\documentclass[3p]{elsarticle} 

\usepackage[utf8]{inputenc}
\usepackage{MasudGroupStyle}

\biboptions{numbers,sort&compress}

\makeatletter
\def\ps@pprintTitle{%
 \let\@oddhead\@empty
 \let\@evenhead\@empty
 \def\@oddfoot{}%
 \let\@evenfoot\@oddfoot}
\makeatother

\begin{document}
\begin{frontmatter}

\title{Error Estimates and Physics Informed Augmentation of Neural Networks for Thermally Coupled Incompressible Navier Stokes Equations}

\author[add1]{Shoaib Goraya\fnref{firstfoot}}
\author[add2]{Nahil Sobh\fnref{secondfoot}}
\author[add1]{Arif Masud\corref{cor1}\fnref{thirdfoot}}
\ead{amasud@illinois.edu}

\fntext[firstfoot]{Graduate Research Assistant}
\fntext[secondfoot]{Affiliate Professor}
\fntext[thirdfoot]{John and Eileen Blumenschein Professor of Mechanics and Computations}

\cortext[cor1]{Corresponding Author}

\address[add1]{Department of Civil and Environmental Engineering, University of Illinois at Urbana Champaign, IL 61801 USA}
\address[add2]{Center of Artificial Intelligence Innovation at National Center of Supercomputing Applications, University of Illinois at Urbana Champaign, IL 61801 USA}

\begin{abstract}

    Physics Informed Neural Networks (PINNs) are shown to be a promising method for the approximation of Partial Differential Equations (PDEs). PINNs approximate the PDE solution by minimizing physics-based loss functions over a given domain. Despite substantial progress in the application of PINNs to a range of problem classes, investigation of error estimation and convergence properties of PINNs, which is important for establishing the rationale behind their good empirical performance, has been lacking. This paper presents convergence analysis and error estimates of PINNs for a multi-physics problem of thermally coupled incompressible Navier-Stokes equations. Through a model problem of Beltrami flow it is shown that a small training error implies a small generalization error. \textit{Posteriori} convergence rates of total error with respect to the training residual and collocation points are presented. This is of practical significance in determining appropriate number of training parameters and training residual thresholds to get good PINNs prediction of thermally coupled steady state laminar flows. These convergence rates are then generalized to different spatial geometries as well as to different flow parameters that lie in the laminar regime. A pressure stabilization term in the form of pressure Poisson equation is added to the PDE residuals for PINNs. This physics informed augmentation is shown to improve accuracy of the pressure field by an order of magnitude as compared to the case without augmentation. Results from PINNs are compared to the ones obtained from stabilized finite element method and good properties of PINNs are highlighted.
    
    \keyword{Machine Learning, PINNs, Thermally Coupled Flows, Physics Informed Augmentation, Navier-Stokes Equations, Neural Networks, Error Estimates, Convergence}
    \endkeyword
\end{abstract}

\end{frontmatter}
\vspace{-0.5cm}

\section{Introduction}\label{Sec: Intro}
    Thermally coupled Navier-Stokes equations serve as a model for a range of flows in engineering and natural sciences. The complexity of flow fields invariably requires high fidelity computational methods that come with their associated cost of computation. Traditional numerical methods come with the challenges of elaborate computer codes for involved mathematical formulations, prohibitive computational cost for resolving all the physical scales, and uncertainty in the numerical values of the parameters in a problem under consideration. In these situations, physics-based data driven techniques can play an important role in filling the knowledge gap between the physical phenomenon and the modeled physics.
    
    Weighted Residual Methods (WRM), when introduced a century ago, provided approximate semi-analytical solution to Partial Differential Equations (PDEs) \cite{ames1967nonlinear,finlayson2013method}. Their popularity stemmed from the fact that approximate solutions (trial functions) were generated by forming a linear combination of known functions and unknown coefficients which could vary over the domain of interest. The number of components in a linear combination were dictated by the size of the resulting system of equations of the unknown coefficients. Since the introduction of WRM, exploring with various function spaces and weights has led to what is known today as the finite difference and the finite element methods as well as meshless and mesh free methods \cite{fletcher1988computational}. Neural networks with unknown functions (activations) and unknown coefficients (weights and biases) have also been proposed to solve linear and nonlinear partial differential equations via the WRM. WRM-based neural network trial function approximations of PDE solution first appeared in 1990's in the works of \cite{dissanayake1994neural,lagaris1998artificial}. The use of Physics Informed Neural Networks (PINNs) to approximate partial differential equations was popularized by the work of \cite{raissi2019deep,raissi2019physics} who leveraged PDE based loss functions along with TensorFlow’s deep learning library to automate the differentiation and subsequently solve for the unknown weights and functions.
    
    In Deep Neural Networks (DNN), it is assumed that an underlying function maps inputs to outputs in the target domain. The functional map is determined by minimizing the loss function via optimizing network's parameters during the process called training. For instance, given a loss function and a class of neural networks \(\mathcal{N}^m\) with \(m\) parameters, the goal is to find a mapping \(\hat{\mathcal{N}}_S({\hat{\bf{W}}},\hat{\bf{b}}): \bfX \to \bfY\) between inputs \(\bfX=\bigcup\limits_{i=1}^{S} \{ x_{i},y_{i}\} \) and outputs \(\bfY=\bigcup\limits_{i=1}^{S} \{ \bfu_{i},p_{i}\} \) such that it minimizes the loss function over a training dataset of \(S\) collocation points. The mapping \(\hat{\mathcal{N}}_S{(\hat{\bf{W}},\hat{\bf{b}})} \in \mathcal{N}^m \) represents trained DNN approximation function with optimized weights \(\hat{\bf{W}}\) and biases \(\hat{\bf{b}}\). \(\hat{\mathcal{N}}_S\) can also be understood as a global approximation function used in the spectral methods. The structure of such DNN consists of an input-layer followed by one or more internal layers ending with an output layer where each layer consists of neurons. For the input layer the neurons are typically the spatio-temporal locations of the domain points. For example, the input layer for a two-dimensional spatial domain \((x,y)\) will consist of two inputs or neurons \({\bf{X}}=\{x,y\}\). The internal layers (also known as hidden layers) can have a variable number of neurons. These neurons represent the unknown functions used in building the WRM trial approximation while the output layer neurons contain the trial function values. The values of any layer \(i\) in the network can be represented in an algebraic form as \({{\bf{z}}_l} = \sigma \left( {{{\bf{W}}_l}{{\bf{z}}_{l - 1}} + {{\bf{b}}_l}} \right) \; {\rm{for}} \; 1 \le l \le L\). The values in the output layer are computed as \({\bf{Y}} = \sigma \left( {{{\bf{W}}_{L + 1}}{{\bf{z}}_L} + {{\bf{b}}_{L + 1}}} \right)\) where \(L\) is the total number of hidden layers, \(\bf{W}\) is the tensor of weights, \(\bf{b}\) is the bias vector, \(\bf{Y}\) is the output vector, and \(\sigma\) is the activation function that acts as a transfer function to propagate data in the network.
    
    The suitability of the DNNs for approximating PDEs arises from their universal approximation properties as well as automatic differentiation \cite{baydin2018automatic}. The universal approximation theorem states that any compactly supported continuous function can be approximated by a neural network generated function to an arbitrary precision \cite{hornik1989multilayer}. It was shown in \cite{pinkus1999approximation} that a single layer neural network of sufficiently large width can uniformly approximate a function and its derivative. Furthermore, the availability of open-source deep learning libraries such as TensorFlow \cite{abadi2016ghe} and PyTorch \cite{paszke2017automatic} makes the implementation straightforward. In the context of solving general PDEs using DNN, let us consider a second order PDE in an open-bounded domain \(\Omega  \in {\mathbb{R}^{{n_{sd}}}}\), where \(n_{sd}\) is the number of spatial dimensions, with non-homogeneous Dirichlet boundary conditions applied at its boundary \(\Gamma  = \partial \Omega\).
    \begin{equation}\label{eq:3}
    \mathcal{L}u\left( x \right) = f \quad \rm{on} \;  \Omega\\
    \end{equation}
    \begin{equation}\label{eq:4}
    u\left( x \right) = g \quad \rm{on} \; \Gamma
    \end{equation}
    where \(\mathcal{L}\) is the differential operator, \(u\) is the unknown solution field, and \(f\) is the source term. The loss function is comprised of the governing equation which is based on some physical balance laws.
    \begin{equation}\label{eq:5}
    {\rm{Loss}} = {\left\| \mathcal{R}^D \right\|}_\Omega ^2 + {\left\| \mathcal{R}^B \right\|}_\Gamma ^2 = \left\| { \mathcal{L}\hat u\left( {x;\bf{W,b}} \right) - f} \right\|_\Omega ^2 + \left\| { \hat u\left( {x;\bf{W,b}} \right) - g} \right\|_\Gamma ^2
    \end{equation}
    where \(\hat u\left( {x;\bf{W,b}} \right)\) is the DNN approximation to the exact solution \(\tilde{u}\left( x \right)\), \(\mathcal{R}^D\) is the domain residual, and \(\mathcal{R}^B\) is the boundary residual. Eq. \ref{eq:5} yields a constrained optimization problem in which the residual or the loss is minimized using optimization strategies such as the stochastic gradient descent methods during the training process:
    \begin{equation}\label{eq:6}
     \left\{{{\hat{\bf{W}},\hat{\bf{b}}}} \right\} = \argmin_{{\bf{W}},{\bfb}} \mathcal{\{R\}}
    \end{equation}
    where \(\hat{\bf{W}} {\rm{and}}\: {\hat{\bf{b}}}\) are optimized weights and biases, respectively, and are used by the DNN to approximate the PDE solution.
    
\begin{figure}[H]
    \begin{center}
    \includegraphics[width=.6\linewidth]{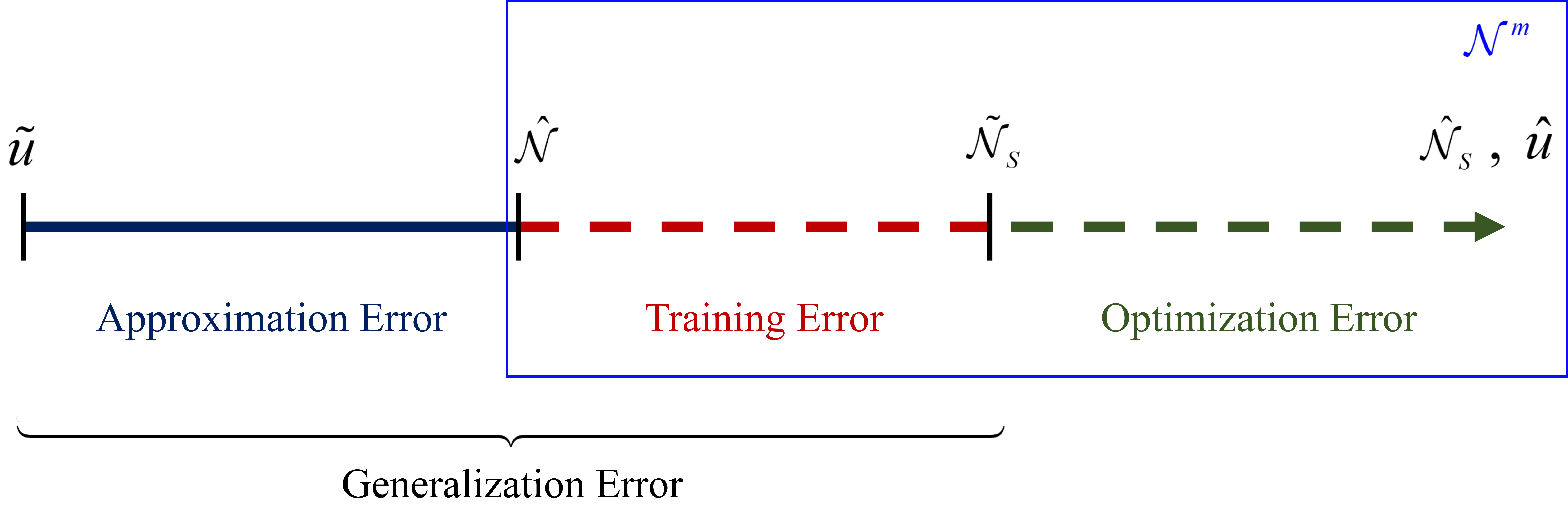}
    \end{center}
    \caption{Schematic diagram of error decomposition in neural network approximated solutions.}
    \label{fig: err_decomp}
    \end{figure}

    The error in supervised learning through neural networks can be decomposed into three components: (i) optimization error \(e_O\), (ii) approximation error \(e_A\), and (iii) training/ estimation error \(e_T\). Fig. \ref{fig: err_decomp} shows a schematic diagram of total error decomposition. Optimization error arises from the minimization problem \ref{eq:4} and is not well understood because the objective function \ref{eq:3} is highly non-convex with multiple local minima. Although several gradient descent algorithms such as \textit{Adam} and \textit{BFGS} have been successfully employed to obtain solutions with enough accuracy, establishing the convergence of these algorithms to a unique global minimum remains an open problem \cite{shin2020convergence}. Approximation error describes the discrepancy between the exact solution and the neural network mapping function on a given network architecture. It is relatively well understood in the light of universal approximation theorem \cite{hornik1989multilayer} as discussed earlier. Training or estimation error \(e_T\) arises when the network is trained on a finite dataset \(S\) to get a mapping \(\hat{\mathcal{N}}_S\) on the target domain. The generalization error \(e_G\) is the combination of approximation and training error and defines the accuracy of the neural network predicted solution. In PDE problems, the generalization error is the distance between a global minimizer of the loss and the exact solution to the PDE. The total error of a general PDE \ref{eq:3} over some norm \( \left\| . \right \| \) is given as:
    
    \begin{equation}\label{eq:6.1}
    {e} = {\left\| {\tilde u} - {\hat u}  \right\|}
    \end{equation}
    where \({\tilde u}\) and \({\hat u}\) are the exact and the neural network predicted solutions of the PDE, respectively. The generalization and training errors are:
    
    \begin{equation}\label{eq:6.2}
    {e_G} = {\left\| \mathcal{R}^D \right\|}_\Omega ^2 + {\left\| \mathcal{R}^B \right\|}_\Gamma ^2
    \end{equation}
    
    \begin{equation}\label{eq:6.3}
    {e_T} = \sum_{i=1}^{S_D} {{\left\| \mathcal{R}_i ^D \right\|}^2} + \sum_{i=1}^{S_B} {\left\| \mathcal{R}_i ^B \right\|} ^2 = \sum_{i=1}^{S_D} \left\| { \mathcal{L}\hat u\left( {x_i;\bf{W,b}} \right) - f_i} \right\|^2 + \sum_{i=1}^{S_B} \left\| { \hat u\left( {x_i;\bf{W,b}} \right) - g_i} \right\| ^2
    \end{equation}
    where \(S_D\) and \(S_B\) are the total number of training collocation points in domain and boundary, respectively.
    
    In recent years, PINNs have been used widely in modeling fluid flows, Navier-Stokes equations \cite{gao2021phygeonet, jin2021nsfnets, mao2020physics, raissi2019deep, rao2020physics,sun2020physics}, solving stochastic PDEs \cite{karumuri2020simulator}, flows in porous media \cite{he2020physics, tripathy2018deep, zhu2018bayesian}, and cardiovascular systems \cite{kissas2020machine, sahli2020physics}. Despite their remarkable performance, theoretical and numerical investigations on their convergence is lacking. Theoretical studies are limited to linear elliptic and parabolic PDEs \cite{shin2020convergence,mishra2022estimates}. To the best of our knowledge, there has been no study on the convergence of PINNs for coupled system of non-linear PDEs.
    
    Since the aim of PINNs is to minimize the training error and in turn the generalization error, bounding the generalization and total errors in terms of training error and the number of collocation points is key to establish the convergence of PINNs to the exact solution of a given system of PDEs \cite{de2022error}. The first objective in this paper is to carry out a numerical investigation on the convergence of PINNs for thermally coupled incompressible Navier Stokes Equations and provide a rationale for their empirical performance. This involves establishing that for a small generalization error \(\epsilon\) the total error is also small i.e., \(e < \delta(\epsilon)\) for some \(\delta(\epsilon) < \mathcal{O}(\epsilon)\). Since PINNs minimize training error instead of generalization error, the above assumption holds under the premise that the generalization error decreases proportionately with the the training error.
    
    Another issue PINNs face for incompressible Navier-Stokes equations is the poor pressure prediction for velocity-pressure formulations. In fact the pressure prediction of the Beltrami flow problem in \cite{jin2021nsfnets} has an \(L_2\) error of ~13\% which is an order of magnitude higher than the velocity field. In the context of multi-physics problems, a less accurate pressure field in turn affects the accuracy of the velocity field. The second objective of this paper is to improve the pressure prediction by introducing a physics-based augmentation in the loss function of the coupled system of PDEs.
    
    In this paper, we have presented convergence analysis of PINNs applied to 2D steady-state flow of incompressible Navier-Stokes PDEs coupled with a scalar energy PDE and subjected to Dirichlet boundary conditions. Through a model problem of Beltrami flows for which an analytical solution exists, effect of (i) network size, (ii) collocation points, and (iii) training error on total error of the predicted solution is investigated and posteriori error estimates are obtained. The systemic study establishes the convergence and consistency of PINNs for the multiphysics problems. Moreover, we introduce a physics-informed augmentation to PINNs in the form of a pressure Poisson equation and its effect on the prediction of the pressure field is shown.
    
    The first contribution in this paper is the convergence analysis and error estimates of PINNs for a multiphysics problem of thermally coupled incompressible Navier-Stokes equations. We show that with enough number of network parameters, collocation points, and a tighter training convergence criterion that leads to a small training error, the PINNs predicted solution converges to the exact solution in various Sobolev norms. Under these conditions, a small training residual implies a small total error. We present \textit{posteriori} convergence rates of total error w.r.t training residual and collocation points. This is of practical importance in determining adequate number of training parameters and residual thresholds to get a good PINNs prediction of coupled thermofluidics for steady state laminar flows. We also show that these convergence rates can be generalized for different geometries and higher Reynolds number flows in the laminar regime.
    
    The second major contribution of the paper deals with physics-informed augmentation of PINNs to get better prediction of pressure field in the coupled system. Borrowing ideas from stabilized methods \cite{masud2004multiscale}, we introduce a pressure stabilization type term in the form of pressure Poisson equation into the PDE residuals for PINNs. This term is  shown to improve the accuracy of the pressure field by an order of magnitude as compared to the case without augmentation. More sophisticated approaches to stabilized methods for Navier-Stokes equations are presented in \cite{masud2006multiscale, masud2009variational,masud2018variationally}, and to turbulent flows in \cite{calderer2013residual, zhu2021residual}. A class of variationally derived stabilized methods for thermofluidic systems are presented in \cite{zhu2021residual, zhu2022nano}. In the present work we have opted to employ a simple approach to keep the focus of the paper on the ML issues. Moreover, we do not employ any penalty parameter for the boundary residuals in the total residual as it is often a tedious trial and error procedure to find a suitable value for a particular problem. This renders the proposed loss function free of any user-defined penalty parameters and is shown to perform robustly in finding a minimzer for the coupled system of PDEs.
    
    An outline of the paper is as follows; Section \ref{Sec: PINNforThermNSE} presents a physics informed neural network in the context of thermally coupled steady state incompressible Navier Stokes equations. We propose a physics-informed augmentation to the standard PINNs method in Section \ref{Sec: PhyInfAugmPINN}. In Section \ref{Sec: NumericExp}, the solution to a subclass of Beltrami flows with a known analytical solution is predicted using the developed framework and an error analysis is carried out. The effect of physics-informed augmentation, network architecture, global training loss, and collocation points on the prediction error is analyzed. The error estimates are shown to have been generalized at different geometries and higher Reynolds's numbers in Section \ref{Sec: Generalization}. Moreover, the PINNs predicted solution is compared with the solution from a stabilized finite element method \cite{zhu2022nano} in Section \ref{SubSec: CompFEMandPINN}. Finally, concluding remarks are presented in Section \ref{Sec: Conclusion}.
\section{PINNs for Thermally Coupled Steady State Navier Stokes Equations} \label{Sec: PINNforThermNSE}
    In this section, we develop PINNs model for solving thermally coupled steady state incompressible Navier-Stokes equations. In thermally coupled flows, the spatially varying active temperature field \(T\left( x \right)\) results in a variable density field \(\rho \left( x \right)\) in the domain. In the context of steady-state flow of incompressible Newtonian fluids, the effect of density variation is accounted for by introducing a buoyancy force in the momentum balance equation, while the continuity equation appears as an incompressibility condition. In addition, a background thermal gradient is applied via thermal boundary condition, and it serves as the driving mechanism for the thermal phase. The governing system of equations for density-stratified flows known as the Boussinesq equations \cite{doering1995applied} defined in an open bounded domain \(\Omega  \in {\mathbb{R}^{{n_{sd}}}}\), can be written as follows.
    \begin{equation}\label{eq:7}
    {\bf{u}} \cdot \nabla {\bfu} + \nabla p - \nabla  \cdot \left( {2\nu {\nabla ^s}{\bfu}} \right) + {\bf{g}}\beta \theta  = {{\bf{f}}_b}
    \end{equation}
    \begin{equation}\label{eq:8}
    \nabla  \cdot {\bf{u}} = 0
    \end{equation}
    \begin{equation}\label{eq:9}
 	{\bf{u}} \cdot \nabla \theta  - \nabla  \cdot \left( {\alpha \nabla \theta } \right) = f
    \end{equation}
    where \(\bfu = \left[ {u,v} \right]\) and \(p\) are the velocity and kinematic pressure fields in two dimensions, respectively; \(\theta  = T - {T_0}\) is the relative temperature (\(T\) is the absolute temperature and \({T_0}\) is the reference temperature); \(\nu \) is the kinematic viscosity of the fluid, \(\beta \) is the thermal expansion coefficient, \({\bf{g}} = \left[ {{g_x},{g_y}} \right]\) is the gravity acceleration vector, \(\alpha \) is the thermal diffusivity; \({{\bf{f}}_b} = \left[ {{f_{bx}},{f_{by}}} \right]\) is the non-gravitational body force, and \(f\) is the heat source/sink. \({\nabla ^s} = {{\left( {\nabla  + {\nabla ^T}} \right)} \mathord{\left/ {\vphantom {{\left( {\nabla  + {\nabla ^T}} \right)} 2}} \right. \kern-\nulldelimiterspace} 2}\) is the symmetric gradient operator. Equation \eqref{eq:7} is the momentum balance equation with a buoyancy term accounting for the thermal effects, Eq. \eqref{eq:8} is the continuity equation to enforce the incompressibility condition, and Eq. \eqref{eq:9} is the energy conservation in the form of the convection-diffusion of the relative temperature field. The boundary conditions on the domain boundary \(\Gamma \) are:
    \begin{equation}\label{eq:10}
 	{\bf{u}}\left( {\bf{x}} \right) = {\bf{g}}_{M} \quad {\rm{on}} \; \Gamma_g^M
    \end{equation}
    \begin{equation}\label{eq:11}
 	\theta \left( {\bf{x}} \right) = {g_E} \quad \rm{on} \; \Gamma_g^E
    \end{equation}
    \begin{equation}\label{eq:12}
 	 {\bf{\sigma}} \cdot {\bf{n}} = \left( {2\nu {\nabla ^s}{\bf{u}} - p{\bf{I}}} \right) \cdot {\bf{n}} = {{\bf{h}}_M} \quad \rm{on} \; \Gamma_h^M
    \end{equation}
    \begin{equation}\label{eq:13}
 	\phi  \cdot n = \alpha \nabla \theta  \cdot n = {h_E} \quad \rm{on} \; \Gamma_h^E
    \end{equation}
    where \({{\bf{g}}_M}\) and \({g_E}\) are the Dirichlet boundary conditions for the velocity field and the relative temperature field, while \({{\bf{h}}_M}\) and \({h_E}\) are the Neumann boundary conditions for the total stress \(\sigma\) and heat flux \(\phi\), respectively. \({\bf{n}}\) is the unit outward normal vector at the boundary. Moreover, these boundaries satisfy the following conditions: \(\Gamma _g^M \cap \Gamma _h^M = \emptyset \), \(\Gamma _g^M \cup \Gamma _h^M = \Gamma \), \(\Gamma _g^E \cap \Gamma _h^E = \emptyset \) and \(\Gamma _g^E \cup \Gamma _h^E = \Gamma \).
    
\begin{figure}[h!]
    \begin{center}
    \includegraphics[width=\textwidth]{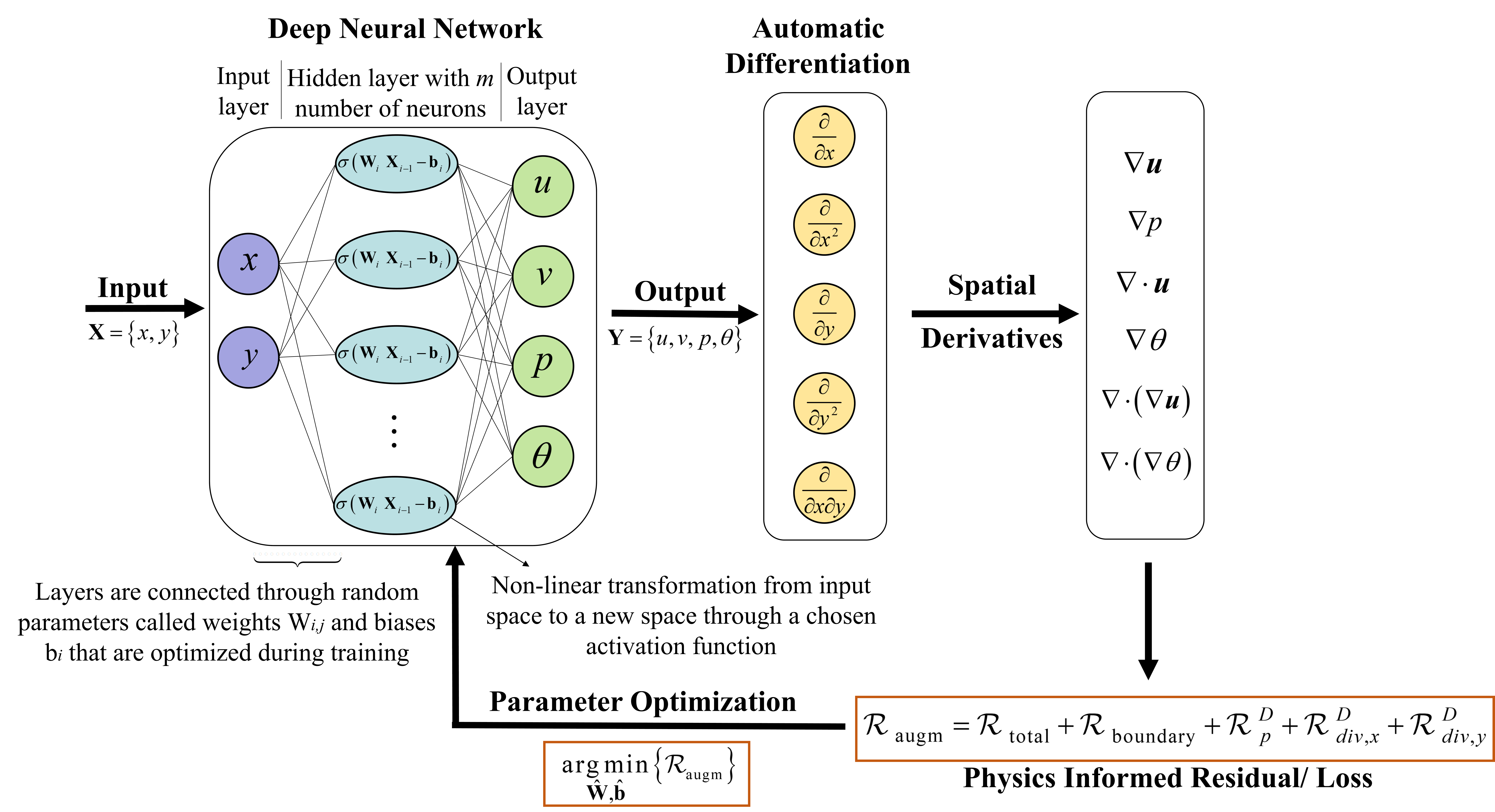}
    \end{center}
    \caption{Schematic diagram of a Physics Informed Neural Networks (PINNs).}
    \label{fig: pinn_schem}
    \end{figure}

	For PINNs in two-dimensions, the inputs to the network are the spatial coordinates \({\bf{X}} = \left\{ {x,y} \right\}\) and the outputs are the four solution fields \({\bf{Y}} = \left\{ {u,v,p,\theta } \right\}\). As opposed to other proposed approaches \cite{Haghighat2021APD} where each solution field is predicted using a separate neural network, we solve the four unknown fields in a coupled fashion using a single network. This approach has advantage in the way that the coupling effect between the mechanical and thermal fields is preserved. The neural network architecture is illustrated in Figure \ref{fig: pinn_schem}.
	
	In this model, only Dirichlet boundary conditions are taken into account. Since PINNs essentially solve an optimization problem, typically the boundary conditions are satisfied by penalizing the residual/ loss function at the boundaries through a Lagrange multiplier \cite{jin2021nsfnets,raissi2019physics,sirignano2018dgm,rao2021physics}, wherein a suitable value of the Lagrange multiplier is chosen through a tedious trial and error procedure. We have not employed any penalty parameter for the boundary residuals, and it is shown through numerical experiments that although the prediction error is mostly accumulated at the boundaries, the boundary residual is less than that of domain during training.
	
	As mentioned in the previous section, a neural network approximates the solution to a PDE by optimized parameters that are obtained by minimizing an objective or loss function. We develop our residual or loss functions using the strong form of the governing PDEs presented in Eq. \eqref{eq:7} – \eqref{eq:9}. In this case, the loss function or the residual is computed using Mean Squared Error (MSE) of the neural network predicted values on the left hand side and the exact values on the right hand side of Eq. \eqref{eq:7} – \eqref{eq:9}. The left hand side of the equations require spatial derivatives and is evaluated using the automatic differentiation \cite{baydin2018automatic}. The domain residual for our PINNs model is given as,
	\begin{equation}\label{eq:14}
 	{\mathcal{R}_{{\rm{domain}}}} = {\mathcal{R}}_u^D + {\mathcal{R}}_v^D + {\mathcal{R}}_{div}^D + {\mathcal{R}}_\theta ^D
    \end{equation}
    where,
    \begin{equation*}
    \begin{split}
    \mathcal{R}_u^D &= {\left\| {\left( {\bfu \cdot \nabla u + \frac{{\partial p}}{{\partial x}} - \nu \Delta u + {g_x}\beta \theta } \right) - {{\left( {{f_{bx}}} \right)}}} \right\|^2},\\
    \mathcal{R}_v^D &= {\left\| {\left( {\bfu \cdot \nabla v + \frac{{\partial p}}{{\partial y}} - \nu \Delta v + {g_y}\beta \theta } \right) - {{\left( {{f_{by}}} \right)}}} \right\|^2},\\
    \mathcal{R}_{div}^D &= {\left\| \nabla  \cdot {\bfu} \right\|^2},\\
    \mathcal{R}_\theta^D &= {\left\| {\left( {\bfu \nabla \theta  - \alpha \Delta \theta } \right) - {{\left( f \right)}}} \right\|^2}
    \end{split}
    \end{equation*}
    whereas boundary residuals are computed by comparing exact values with the
    \begin{equation}\label{eq:15}
    {\mathcal{R}_{{\rm{boundary}}}} = {\mathcal{R}}_u^B + {\mathcal{R}}_v^B + {\mathcal{R}}_\theta^B
    \end{equation}
    where,
    \begin{equation*}
    \begin{split}
    {\mathcal{R}}_u^B &= {\left\| { u - {{ {{g_{Mx}}}}}} \right\|^2},\\
    {\mathcal{R}}_v^B &= {\left\| { v - {{ {{g_{My}}}}}} \right\|^2},\\
    {\mathcal{R}}_\theta^B &= {\left\| { \theta - {{ {{g_E}} }}} \right\|^2}
    \end{split}
    \end{equation*}
    The total loss or residual is comprised of the residual from both the domain and the boundaries:
    \begin{equation}\label{eq:16}
    {\mathcal{R}_{{\rm{total}}}} = {\mathcal{R}_{{\rm{domain}}}} + {\mathcal{R}_{{\rm{boundary}}}}
    \end{equation}
\section{Physics Informed Augmentation to PINNs}\label{Sec: PhyInfAugmPINN}
    For PINNs applied to incompressible Navier-Stokes equations, satisfying incompressibility constraint is an important issue. Since pressure acts as a Lagrange multiplier that enforces the incompressibility constraint, any error arising due to the weak enforcement of this constraint leads to less accurate prediction of the pressure field \cite{jin2021nsfnets}. In order to get a better prediction of the pressure field, strategies involving training of neural networks via an augmented dataset, enforcement of penalty on the pressure field during the training process, and input/ output perturbation are employed \cite{jin2021nsfnets,simard2003best,yaeger1996effective}. We propose a physics-informed augmentation in which the residuals are constructed based on a pressure Poisson equation and spatial derivatives of the incompressibility constraint. The augmentation is intended to help improve the pressure prediction by satisfying additional physics-based constraints within the domain. The governing equations of the augmentation in the domain \(\Omega  \in {\mathbb{R}^{{n_{sd}}}}\) are given as follows,
    \begin{equation}\label{eq:17}
    {\Delta}p = \nabla  \cdot \left( {{\bf{f}}_{b} - {\bf{u}} \cdot \nabla {\bf{u}} - {\bf{g}} \beta \theta } \right)
    \end{equation}
    \begin{equation}\label{eq:18}
    {\left( {\nabla  \cdot {\bf{u}}} \right)_{,x}} = 0
    \end{equation}
    \begin{equation}\label{eq:19}
    {\left( {\nabla  \cdot {\bf{u}}} \right)_{,y}} = 0
    \end{equation}
    This gives rise to an augmented residual or loss function which is given as,
    \begin{equation}\label{eq:20}
    \mathcal{R}_{\rm{augm}} = \mathcal{R}_p^D + \mathcal{R}_{div,x}^D + \mathcal{R}_{div,y}^D
    \end{equation}
    where,
    \begin{align*}
    {\mathcal{R}}_p^D &= {\left\| {\left( {{\Delta}p} \right) - {{\left( {\nabla  \cdot \left( {{\bf{f}}_{b} - \bfu \cdot \nabla \bfu - {\bf{g}}\beta \theta } \right)} \right)}}} \right\|^2},\\
    {\mathcal{R}}_{div,x}^D &= {\left\| {{{\left( {\nabla  \cdot {\bf{u}}} \right)}_{,x}}} \right\|^2},\\
    {\mathcal{R}}_{div,y}^D &= {\left\| {{{\left( {\nabla  \cdot {\bf{u}}} \right)}_{,y}}} \right\|^2}
    \end{align*}
    The proposed augmentation sits in the total residual given in Eq. \eqref{eq:16} and is shown to significantly improve pressure field prediction in the subsequent section.
    \begin{equation}\label{eq:21}
    {\tilde{\mathcal{R}}_{{\rm{total}}}} = {\mathcal{R}_{{\rm{domain}}}} + {\mathcal{R}_{{\rm{boundary}}}} + {\mathcal{R}_{{\rm{augm}}}}
    \end{equation}
\section{Numerical Experiments: Beltrami Flows}\label{Sec: NumericExp}
\begin{figure}[h!]
    \begin{center}
    \includegraphics[width=.8\linewidth]{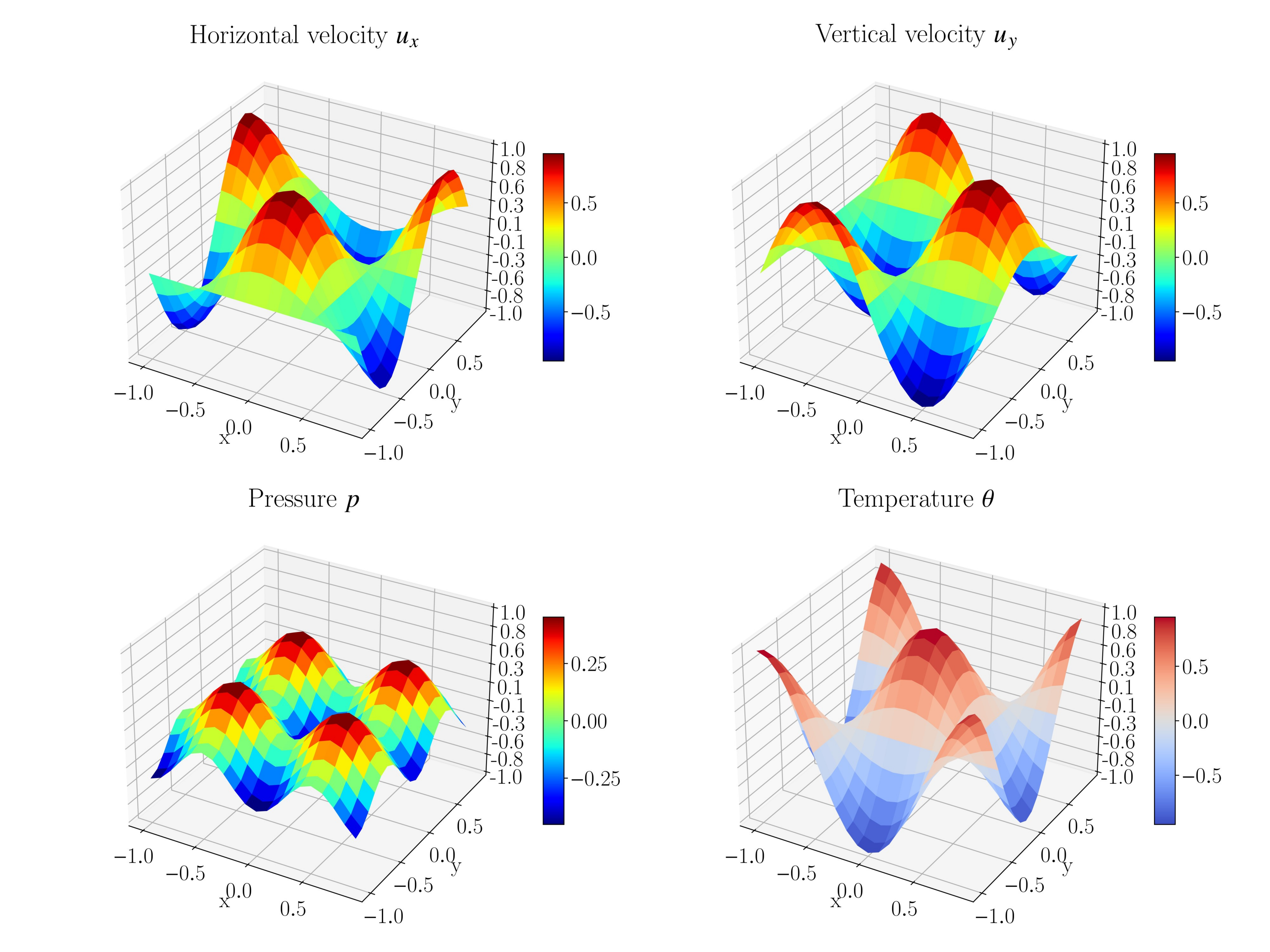}
    \end{center}
    \caption{Visualization of the analytical solution for Beltrami flow: velocity, pressure, and temperature fields.}
    \label{fig: beltrami_soln}
    \end{figure}

    We employ the PINNs model developed in the previous sections to predict solution of a subclass of Beltrami flows \cite{wong2007numerical}. The problem has an analytical solution against which the trained neural networks are tested. This class of problems is used to study buoyancy-induced convection and heat transfer phenomenon and has wide-ranging applications in engineering systems, such as solar collectors, electronic cooling, heat exchangers, and thermal insulating systems etc. \cite{garnier2009integrated}. The analytical expressions for the velocity, pressure and relative temperature fields are as follows,
    \begin{equation}\label{eq:22}
    \bfu \left( {x,y} \right) = {\left[ { - \cos \left( {\pi x} \right)\sin \left( {\pi y} \right),\sin \left( {\pi x} \right)\cos \left( {\pi y} \right)} \right]^T}
    \end{equation}
    \begin{equation}\label{eq:23}
    p\left( {x,y} \right) =  - \frac{1}{4}\left( {\cos \left( {2\pi x} \right) + \cos \left( {2\pi y} \right)} \right)
    \end{equation}
    \begin{equation}\label{eq:24}
    \theta \left( {x,y} \right) = \cos \left( {\pi x} \right)\cos \left( {\pi y} \right)
    \end{equation}
    and the body forces and heat source driving the flow are as follows,
    \begin{equation}\label{eq:25}
    {\bf{f}}_{b}\left( {x,y} \right) = \left[ {\begin{array}{*{20}{c}}{ - 2{\pi ^2}\nu \cos \left( {\pi x} \right)\sin \left( {\pi y} \right) + {g_1}\beta \cos \left( {\pi x} \right)\cos \left( {\pi y} \right)}\\{2{\pi ^2}\nu \sin \left( {\pi x} \right)\cos \left( {\pi y} \right) + {g_2}\beta \cos \left( {\pi x} \right)\cos \left( {\pi y} \right)}\end{array}} \right]
    \end{equation}
    \begin{equation}\label{eq:26}
    f\left( {x,y} \right) = 2{\pi ^2}\alpha \cos \left( {\pi x} \right)\cos \left( {\pi y} \right)
    \end{equation}
    The computational domain is a bi-unit square \(\Omega  = \left[ {-1,1} \right] \times \left[ {-1,1} \right]\) that covers one period of the sinusoidal velocity and pressure solution fields.  Dirichlet boundary conditions \({\bf{g}}\left( {\bf{x}} \right) = {\left. {[\bfu,p,\theta]} \right|_{x  \in {\Gamma}}}\) for the velocity, pressure, and the relative temperature fields are applied at the domain boundaries \(\Gamma  = \partial \Omega \). The non-dimensional material parameters used in our problem are all unity i.e., \(\alpha  = \beta  = \nu  = 1.0\). The exact solution and boundary conditions are shown in Figure \ref{fig: beltrami_soln}.
    
\begin{figure}[H]
        \centering
        \begin{subfigure}[t]{.27\linewidth}
            \centering
            \includegraphics[width=\linewidth]{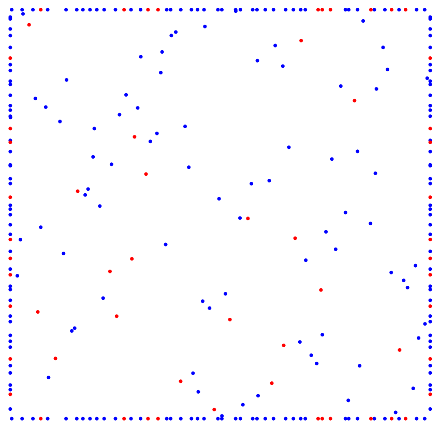}
            \caption{Training and validation points (15\% of the points, shown in red, are used for validation}
        \end{subfigure}
        \hspace{1.5em}%
        \begin{subfigure}[t]{.27\linewidth}
            \centering
            \includegraphics[width=\linewidth]{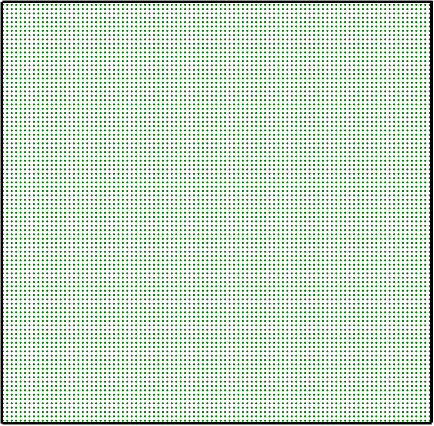}
            \caption{Test points on a \(100 \times 100\) grid.}
        \end{subfigure}
        \caption{Schematic diagrams of data split and random collocation points for training and validation of PINNs.}
        \label{fig: train_test_pts}
    \end{figure}

    For training, random collocation points are generated in the domain as well as on the boundaries using Latin hypercube sampling strategy \cite{mckay2000comparison}. Total collocation points are chosen such that the domain has twice the number of points than on all the four boundaries combined. For convergence analysis, collocation points are hierarchically increased such that the next set of points also contains the previous points. The details of each training dataset is given in Table \ref{table:1}. The trained network is tested on a \(100 \times 100\) uniform grid of collocation points that covers the problem domain and boundaries. A schematic diagram of training and test points is shown in \ref{fig: train_test_pts}.
    
    \begin{table}[H]
    \centering
    \small
    \caption{Description of hierarchical training datasets used in the convergence analysis.}
    \label{table:1}
    \begin{tabular}{@{}cccc@{}}
    \toprule
    Total Points & Domain Points & Total Boundary Points & Points per Boundary \\ \midrule
    12           & 8             & 4                     & 1                   \\
    24           & 16            & 8                     & 2                   \\
    48           & 32            & 16                    & 4                   \\
    96           & 64            & 32                    & 8                   \\
    192          & 128           & 64                    & 16                  \\
    384          & 256           & 128                   & 32                  \\
    768          & 512           & 256                   & 64                  \\
    1536         & 1024          & 512                   & 128                 \\ \bottomrule
    \end{tabular}
    \end{table}
    
\subsection{Convergence Analysis and Error Estimates} \label{SubSec: ErrorEstimates}
    In order to answer the questions put forth in Section \ref{Sec: Intro}, we systematically evaluate the impact of the choice of network architecture, physics-informed augmentation, collocation points, and training error on the total error. This methodology enables us to numerically estimate \textit{a-posterior} convergence rates of the total error that are of practical importance in selecting appropriate PINNs model for thermally coupled Navier-Stokes equations. 
    
    For this purpose, the solution \(\hat u \) from a trained network is predicted on the \(100 \times 100\) grid of testing points and the total error is computed w.r.t analytical solution \(\tilde{u}\) in Sobolev norms as defined below:
    \begin{equation}\label{eq:27}
    \left\| {\tilde{u} - \hat{u}} \right\|_{{W^{k,\infty}}\left( \Omega \right)} = \mathop{{{\max}}}\limits_{{0 \le m \le k}} {| {\tilde{u} - \hat{u}} |}_{{W^{m,\infty}}\left( \Omega \right)}
    \end{equation}
    where,
    \begin{equation*}
    { | {\tilde{u} - \hat{u}} |}_{{W^{m,\infty}}\left( \Omega \right)} = \mathop{{{\max}}}\limits_{|\alpha| = m} {\left\| D^{\alpha} \left( {\tilde{u} - \hat{u}} \right) \right\|}_{{L^{\infty}}\left( \Omega \right)} \quad {\rm{for}} \:\: {m = 0,...,k}
    \end{equation*}
    and \(D^{\alpha}\) is the operator for spatial derivative of order \(\alpha\).
    
    The total error on test points ensures that the error estimates are also applicable to the data not used in the training and without any prior knowledge of the exact solution.
\subsubsection{Effect of Network Architecture on Approximation and Training Errors} \label{SubSec: EffNNArch}
\begin{figure}[H]
    \centering \includegraphics[width=.65\linewidth]{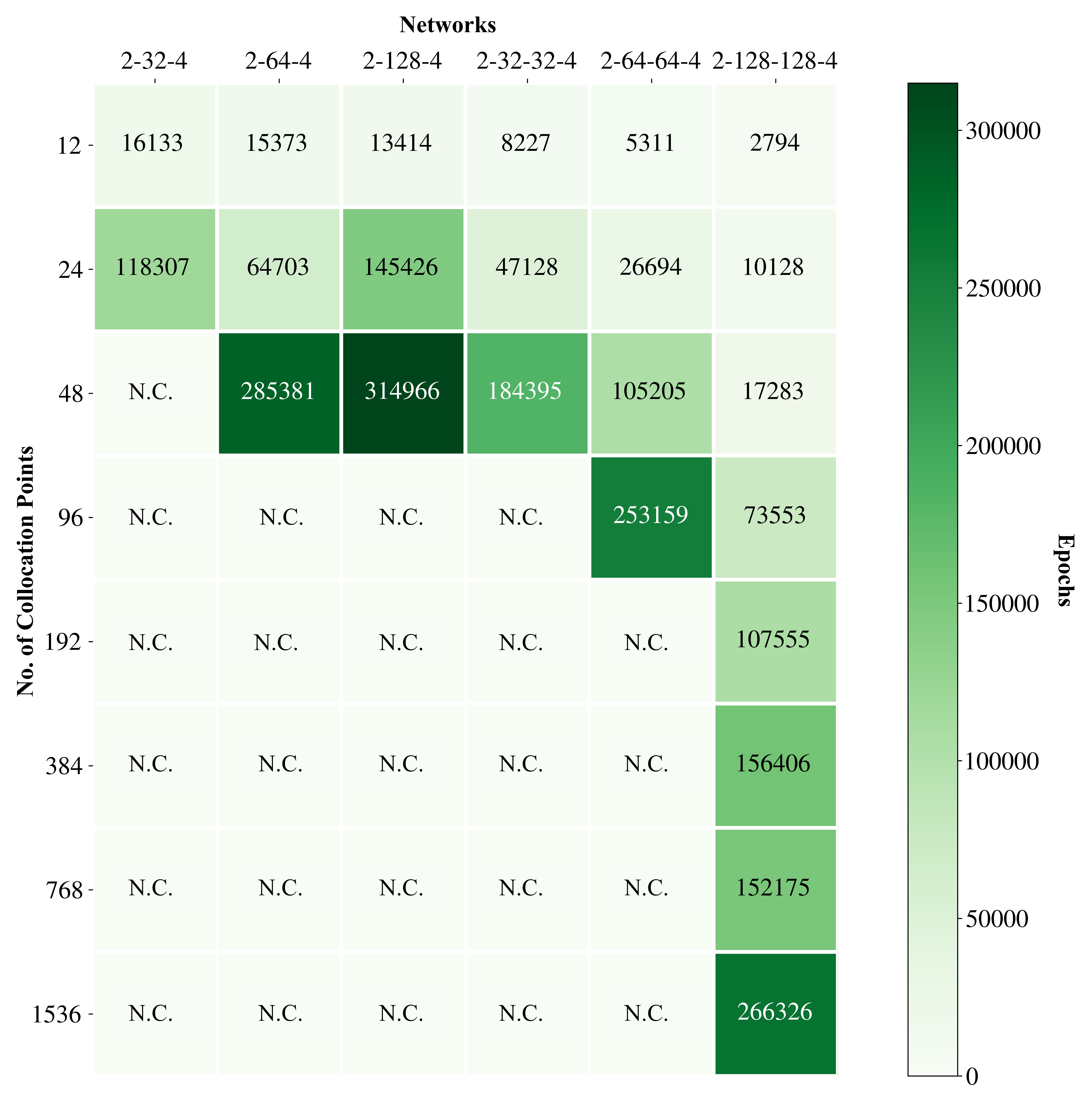}
    \caption{\small Heat map showing the effect of network architecture on the approximation capability of PINNs (N.C. = Network did not converge to the specified training residual of \(10^{-4}\) and exceeded the maximum specified epochs of 350K).}
    \label{fig: heatmap}
    \end{figure}

    Since generalization error comprises of both the approximation and training errors, it is important to have a neural network that is a consistent estimator, and is capable of approximating the given loss function with a mapping that converges to a specified threshold value for the training error. This section investigates the effect of network architecture and trainable parameters on approximation and convergence capabilities of PINNs. Six different network architectures are employed with single and double hidden layers and varying number of neurons in each layer (32, 64, and 128). As we increase the number of hidden layers or neurons, the number of trainable parameters of network also increases as shown in Fig. \ref{fig: heatmap}. The stopping criterion for training of each network is either the training error of \({10^{-4}}\) or the maximum epochs of \(350,000\). With this in place, each network is trained on a range of collocation points. The numbers within each cell indicate the epochs it took for a network to reach the specified training error while N.C. indicates that the network has failed to converge to the training threshold and instead exceeded the maximum number of epochs.
    
    From Fig. \ref{fig: heatmap} it is observed that as we increase the number of collocation points, shallower networks with lesser trainable parameters fail to converge. As we go to deeper and wider networks, network parameters as well as the approximation capability increases, and the network converges faster to the training threshold \cite{raghu2017expressive}. This is consistent with the universal approximation theorem. Based on the study, we chose 2-128-128-4 architecture for subsequent error analysis.
    
    \begin{remark} \label{rem: 1}
    One key take away from this study is that the convergence properties of a given neural network depends on the size of training dataset as well as on the training error. For instance, networks of fixed size have irreducible approximation error irrespective of the size of training dataset (beyond a specific threshold), i.e., increasing the number of collocation points beyond this threshold would not reduce total error and networks are no longer consistent estimators.
    \end{remark}
\subsubsection{Effect of Physics Informed Augmentation on Pressure Prediction} \label{SubSec: EffPhyInfAugm}

\begin{figure}[h!]
    \centering
    \includegraphics[width= .4\linewidth]{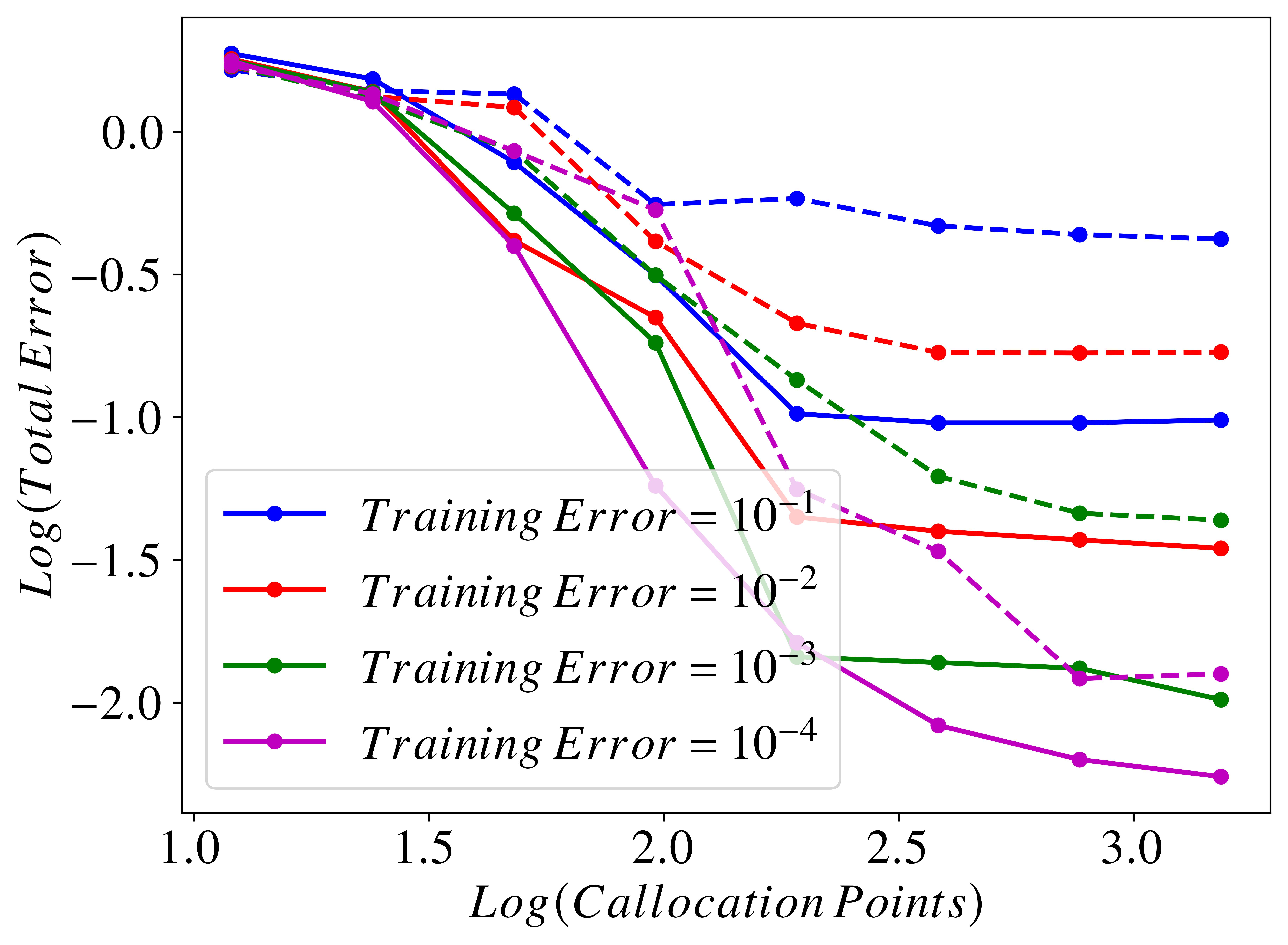}
    \caption{Comparison of error in the pressure field as a function of collocation points for the augmented and the non-augmented cases.}
    \label{fig: pressure_augm_plot}
    \end{figure}

\begin{figure}[h!]
    \centering 
    \includegraphics[width= .63\linewidth]{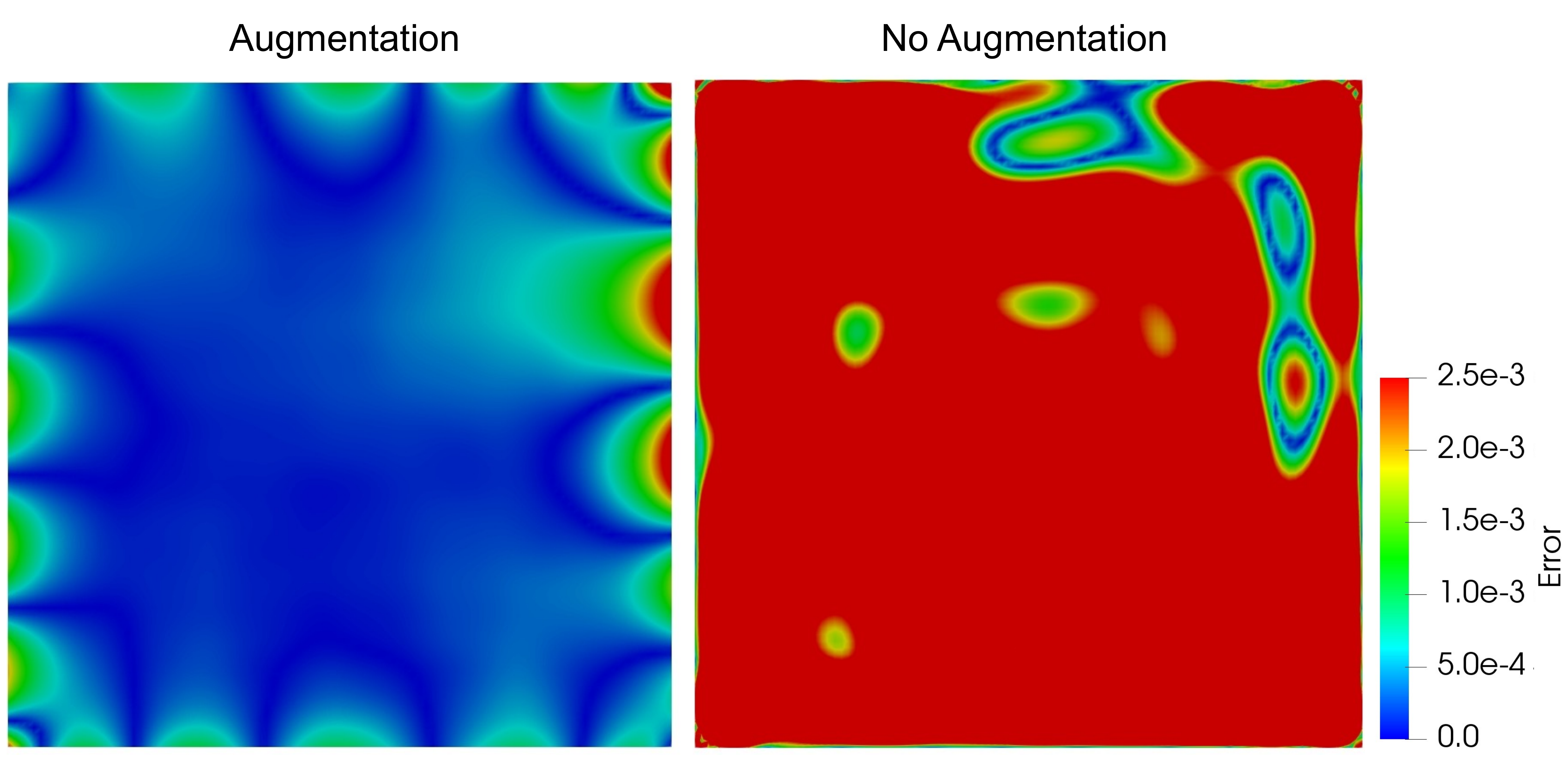}
    \caption{Qualitative plots of pointwise error w.r.t analytical solution of the pressure field for the augmented and the non-augmented cases (training error = \(10^-4\), collocation points = 1536).}
    \label{fig: pressure_augm_qual}
    \end{figure}

    In this section, we investigate the effect of proposed augmentation on PINNs predicted pressure field. Figure \ref{fig: pressure_augm_plot} shows the convergence of total error as a function of collocation points for several training error thresholds. Physics informed augmentation significantly reduces error in the pressure field at every training threshold. This shows that the additional terms in the loss function are significant for a better pressure prediction in the coupled system of PDEs. The observation is further confirmed by the qualitative error plots for pressure field in Figure \ref{fig: pressure_augm_qual}. Without augmentation, the error w.r.t analytical solution is an order of magnitude higher than the case with augmentation.
    
\begin{figure}[H]
    \centering 
    \includegraphics[width= .4\linewidth]{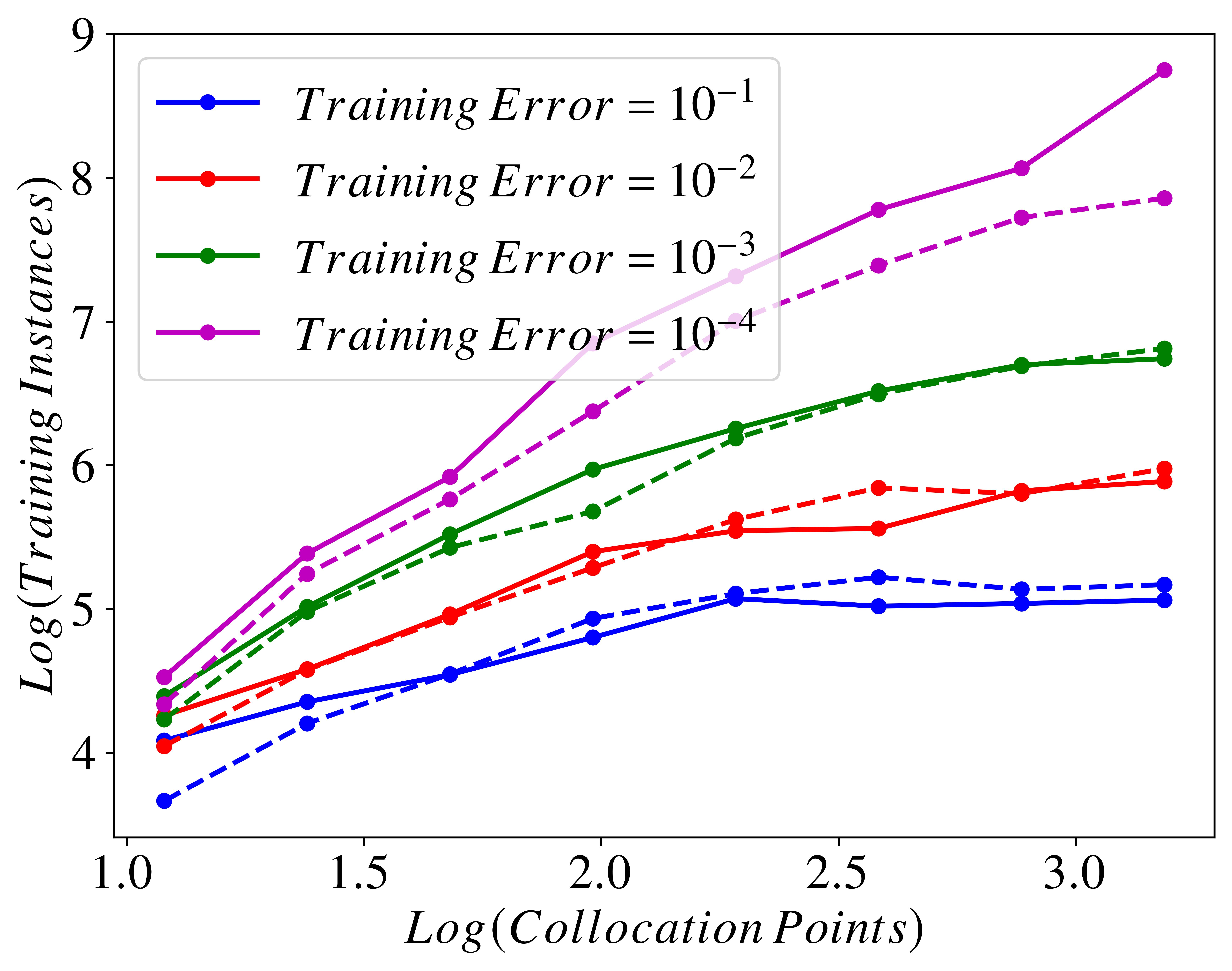}
    \caption{Computational cost for physics informed augmentation.}
    \label{fig: phy_augm_cost}
    \end{figure}

    We also investigate the computational cost of the physics informed augmentation. There is no significant increase in the computational cost with the added augmentation for larger training loss thresholds (\({10^{-1}}\), \({10^{-2}}\), and \({10^{-3}}\)) as shown in Figure \ref{fig: phy_augm_cost}. For the training error of \({10^{-4}}\), computational cost of training PINNs with physics informed augmentation is 2 to 7 times than that of no augmentation. The computational complexity w.r.t the total collocation points for the augmented and the non-augmented case is \(\mathcal{O}\left( {{N^{1.9}}} \right)\) and \(\mathcal{O}\left( {{N^{1.7}}} \right)\),  respectively, where \(N\) is the number of collocation points.
    
\subsubsection{Convergence Analysis of Total Error} \label{SubSec: EffPhyInfAugm}

    As discussed earlier, the underlying idea in PINNs is to minimize the training error over some network parameters. To justify their empirical performance, we have to establish that for a given small training error and in turn the generalization error, the corresponding total error is also small. This section investigates the relation between the training and total error and we estimate the convergence rates that are important in establishing the convergence of PINNs. The network architecture of 2-128-128-4 with augmentation is selected and the total training threshold is varied as \({10^{-1}}\), \({10^{-2}}\), \({10^{-3}}\), and \({10^{-4}}\). For a specified number of collocation points, the network is trained to a given training error and the total error between the exact solution and the PINNs predicted solution is computed over \(W^{0,\infty}\), \(W^{1,\infty}\), and \(W^{2,\infty}\) Sobolev norms. The number of collocation points are hierarchically increased to determine enough number points for optimal convergence.
    
    \begin{figure}[H]
    \centering 
    \includegraphics[width= .75\linewidth]{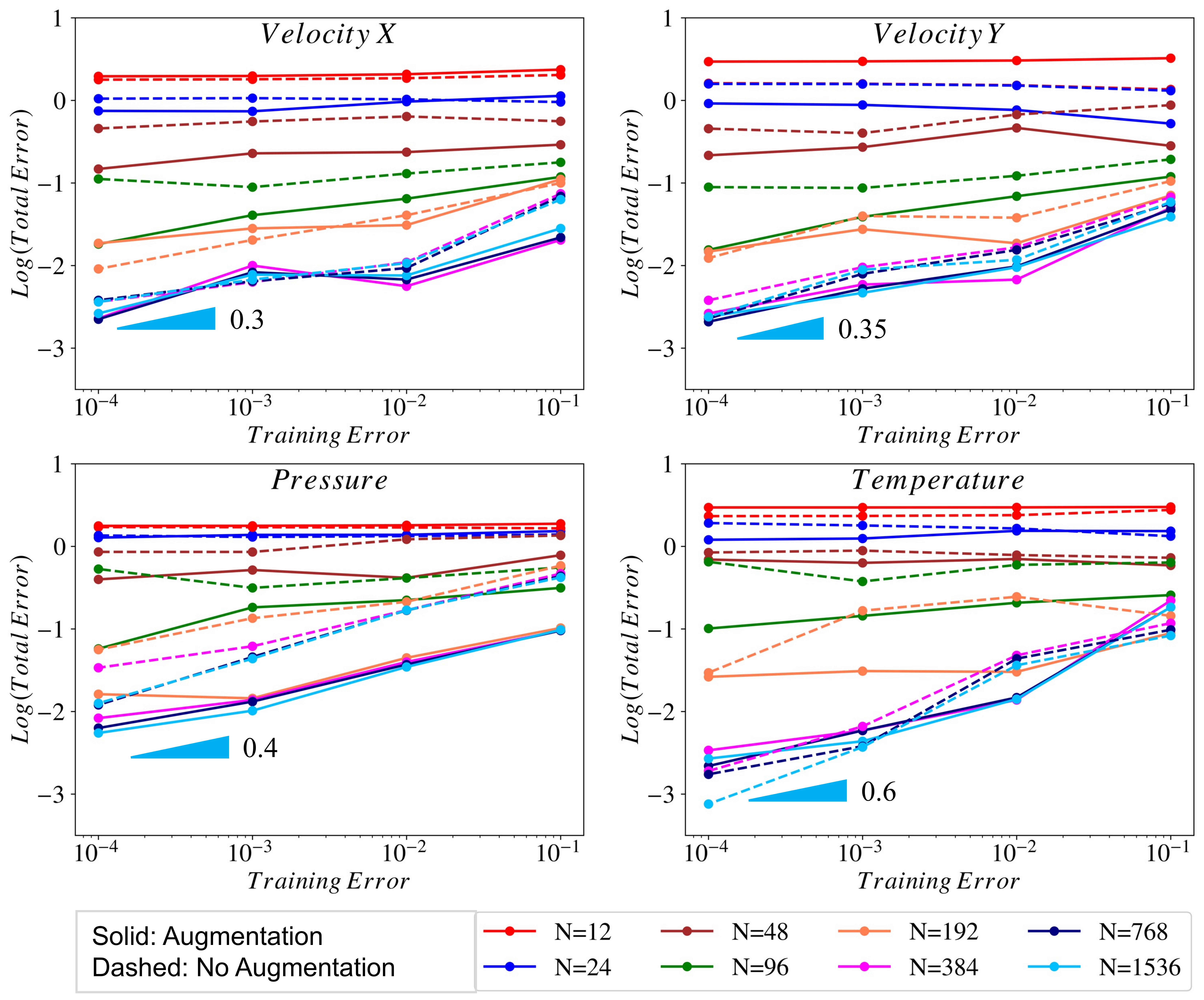}
    \caption{Convergence of \(W^{0,\infty}\) norm of error w.r.t training error.}
    \label{fig: w0_inf_augm}
    \end{figure}
    
    \begin{figure}[h!]
    \centering 
    \includegraphics[width= .7\linewidth]{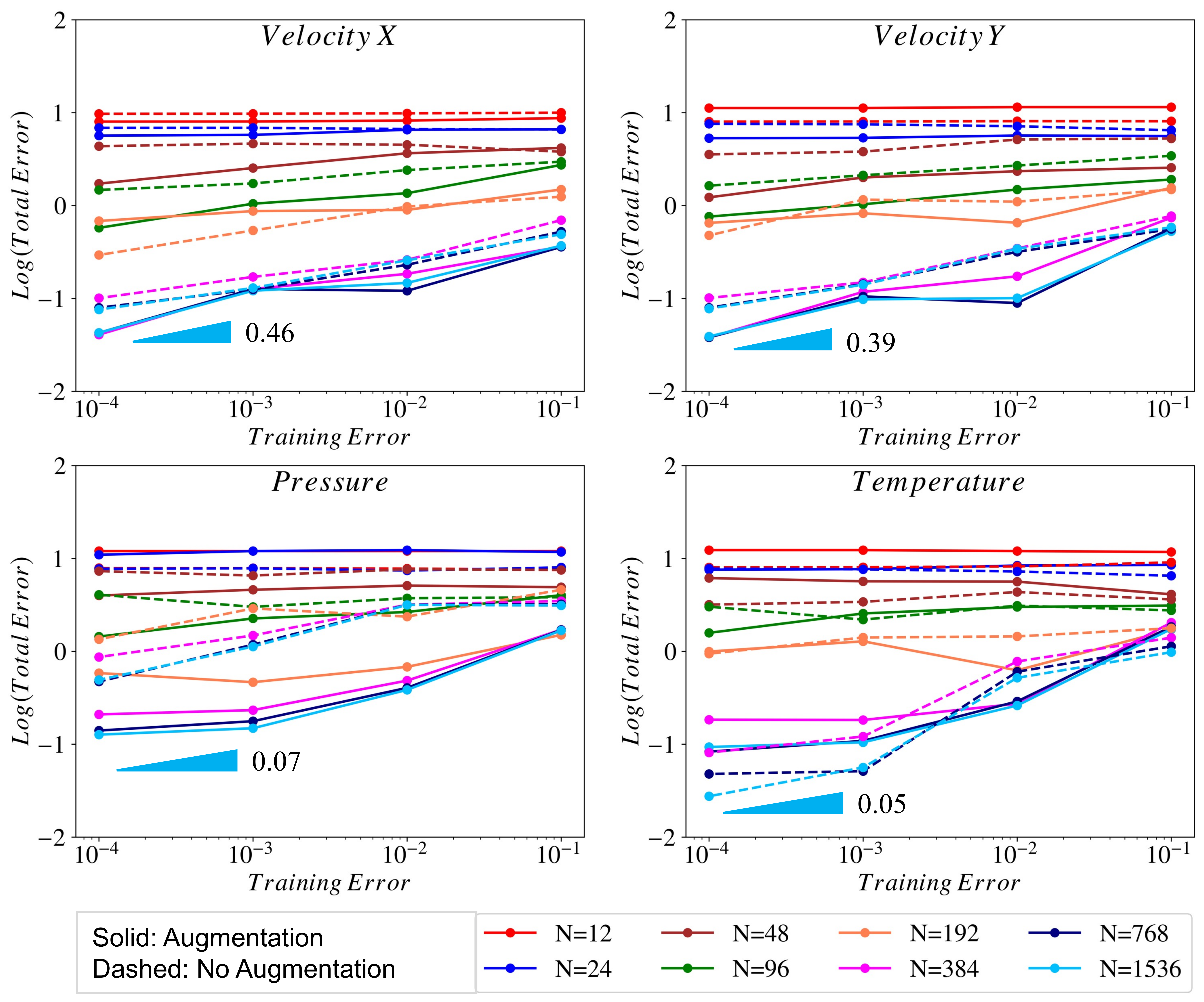}
    \caption{Convergence of \(W^{1,\infty}\) norm of error w.r.t training error.}
    \label{fig: w1_inf_augm}
    \end{figure}
    
\begin{figure}[H]
    \centering 
    \includegraphics[width= .7\linewidth]{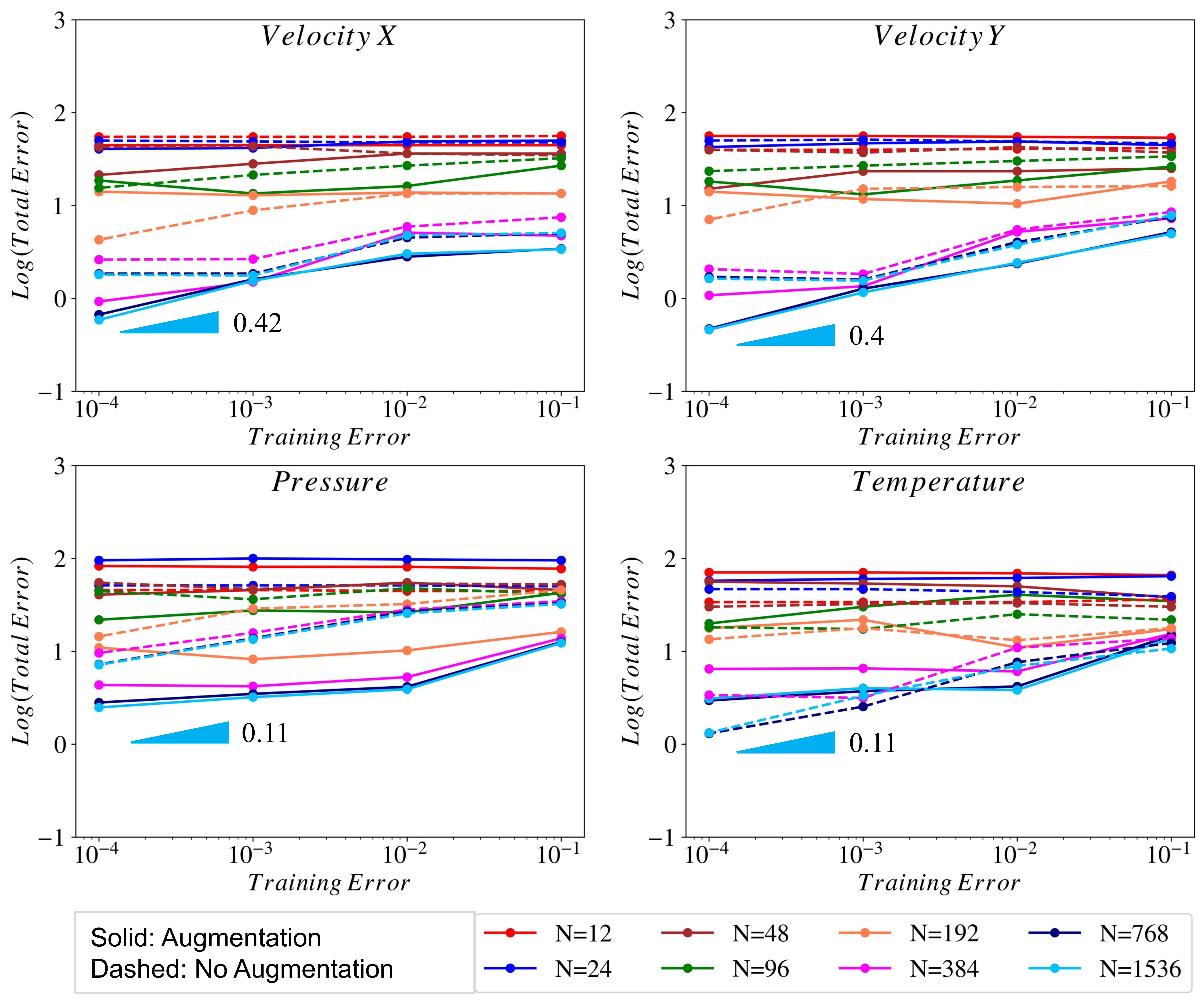}
    \caption{Convergence of \(W^{2,\infty}\) norm of error w.r.t training error.}
    \label{fig: w2_inf_augm}
    \end{figure}

\begin{figure}[H]
        \centering
        \begin{subfigure}[t]{.31\linewidth}
            \centering
            \includegraphics[width=\linewidth]{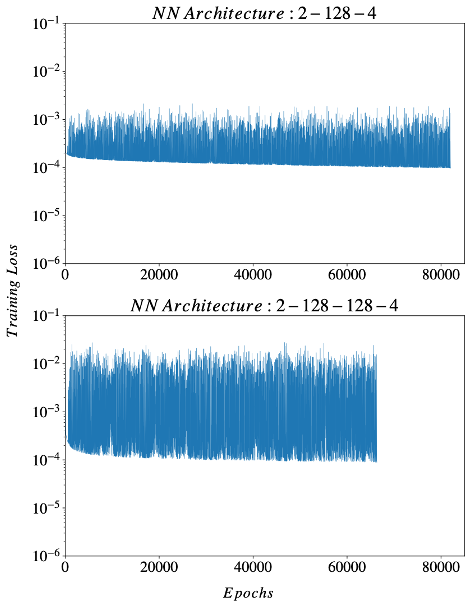}
            \caption{Domain error}
        \end{subfigure}
        \hspace{0.2em}%
        \begin{subfigure}[t]{.31\linewidth}
            \centering
            \includegraphics[width=\linewidth]{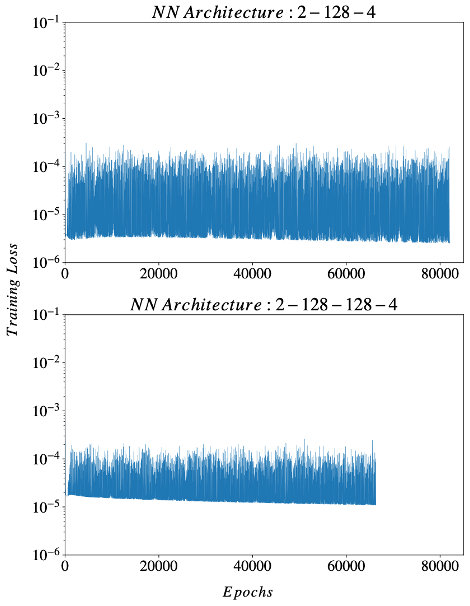}
            \caption{Boundary error}
        \end{subfigure}
        \hspace{0.2em}%
        \begin{subfigure}[t]{.31\linewidth}
            \centering
            \includegraphics[width=\linewidth]{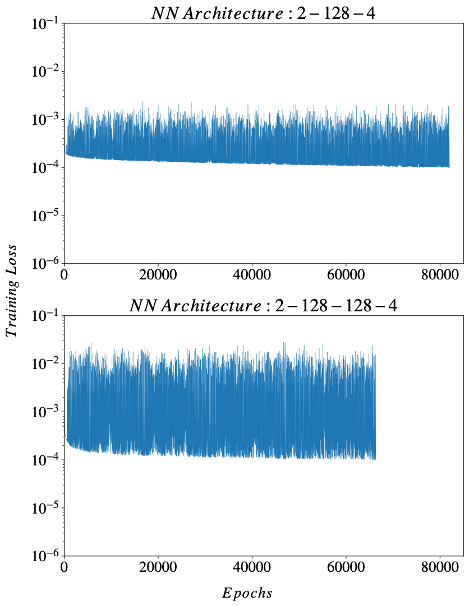}
            \caption{Total error}
        \end{subfigure}
        \caption{Decomposition of training error/ loss for single and double layer neural networks.}
        \label{fig: training_residual_decomp}
    \end{figure}

\begin{figure}[h!]
    \centering 
    \includegraphics[width= .4\linewidth]{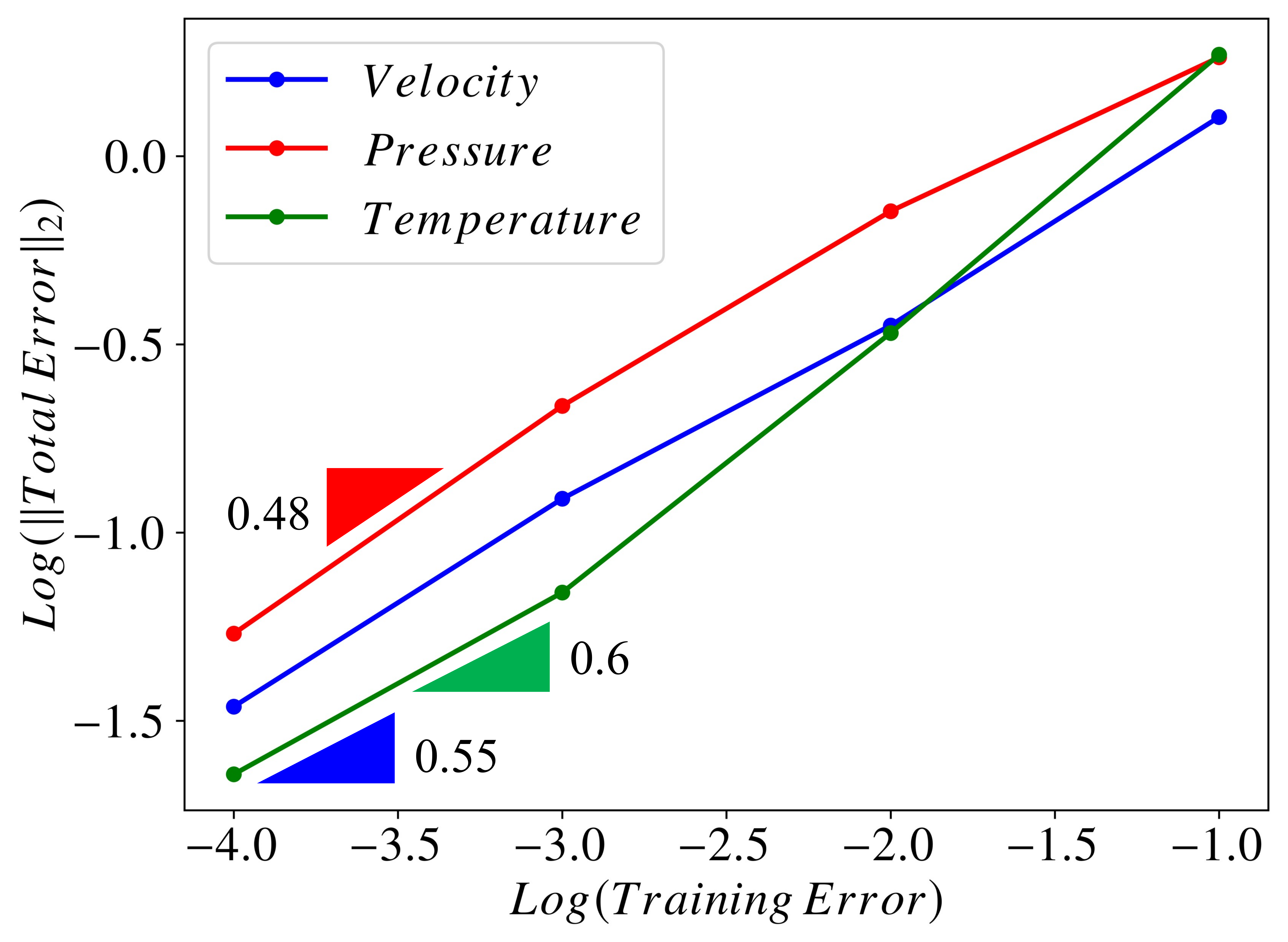}
    \caption{Scaling of total error with training error.}
    \label{fig: l2_error_conv}
    \end{figure}

    It is observed that the number of collocation points are important in getting optimal convergence of total error. All the plots in Figs. \ref{fig: w0_inf_augm}-\ref{fig: w2_conv} show that after 200 points, the convergence rates becomes optimal. For lower number of collocation points, there is no significant convergence. We estimate the convergence rates for the solution fields, first derivatives, and second derivatives in Figs. \ref{fig: w0_inf_augm}-\ref{fig: w2_conv}, respectively. In general the convergence rates improve as we go to higher order derivatives. For 1536 collocation points, the approximate convergence rates of total error for velocity, pressure, and temperature fields are \(\mathcal{O}\left( {{{e_T}^{1/3}}} \right)\), \(\mathcal{O}\left( {{{e_T}^{1/2}}} \right)\), and \(\mathcal{O}\left( {{{e_T}^{1/2}}} \right)\), respectively.
    
     We also analyze the training error in the domain as well as on the boundaries for a single and double layer network in Figure \ref{fig: training_residual_decomp}. The threshold value is set at \({10^{-4}}\). It is seen that the boundary error is an order magnitude less than the domain error which indicates that the total training error is not limited by the weak enforcement of the boundary conditions in this case. Moreover, it is observed that deeper network shows greater oscillations during the training. Since deeper networks have larger number of parameters, these oscillations might be due to redundancy in the parameterisation of the function space in a neighborhood of a local minimum \cite{saarinen1993ill}.
     
    Recently, theoretical error estimates for incompressible Navier-Stokes equations under periodic boundary conditions were presented in \cite{de2022error}. Although the PDE system and the boundary conditions do not exactly match the present case, the error bounds give scaling between total error and training error as follows.
    \begin{equation}\label{eq:29}
    {\int\limits_\Omega  {\left\| {\tilde u  - \hat u} \right\|_2}}\;dx \: \lessapprox \: \mathcal{O} {\left({\sqrt{{e _T}}}\: \right)}
    \end{equation}
    Figure \ref{fig: l2_error_conv} shows that total error scales approximately as the square root of the training error for all the solution fields which validates our numerical findings.
    
\subsubsection{Consistency of PINNs} \label{SubSec: EffCollocPts}

\begin{figure}[H]
    \centering 
    \includegraphics[width= .7\linewidth]{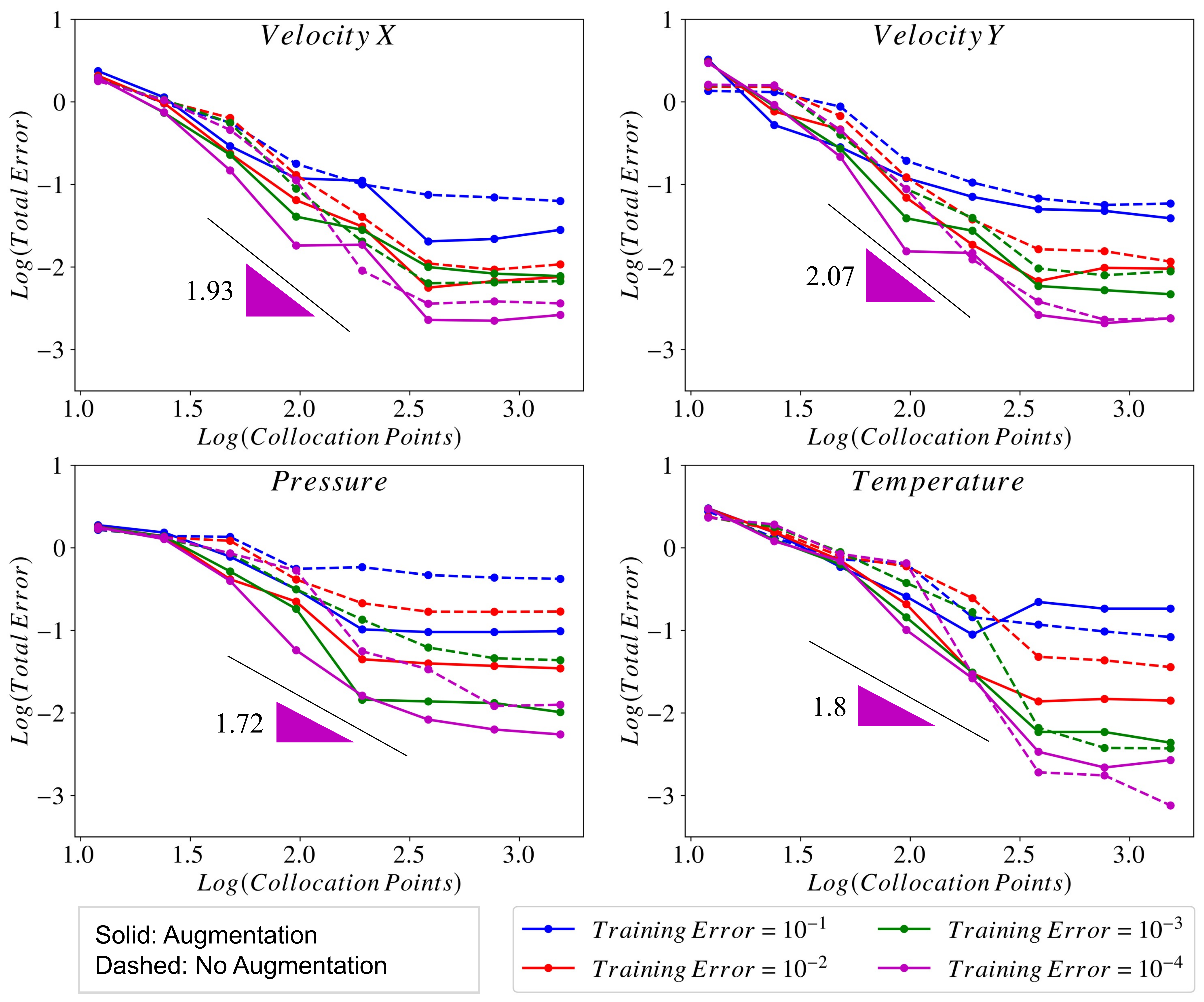}
    \caption{Convergence of \(W^{0,\infty}\) norm of error w.r.t training collocation points.}
    \label{fig: w0_conv}
    \end{figure}

\begin{figure}[H]
    \centering 
    \includegraphics[width= .7\linewidth]{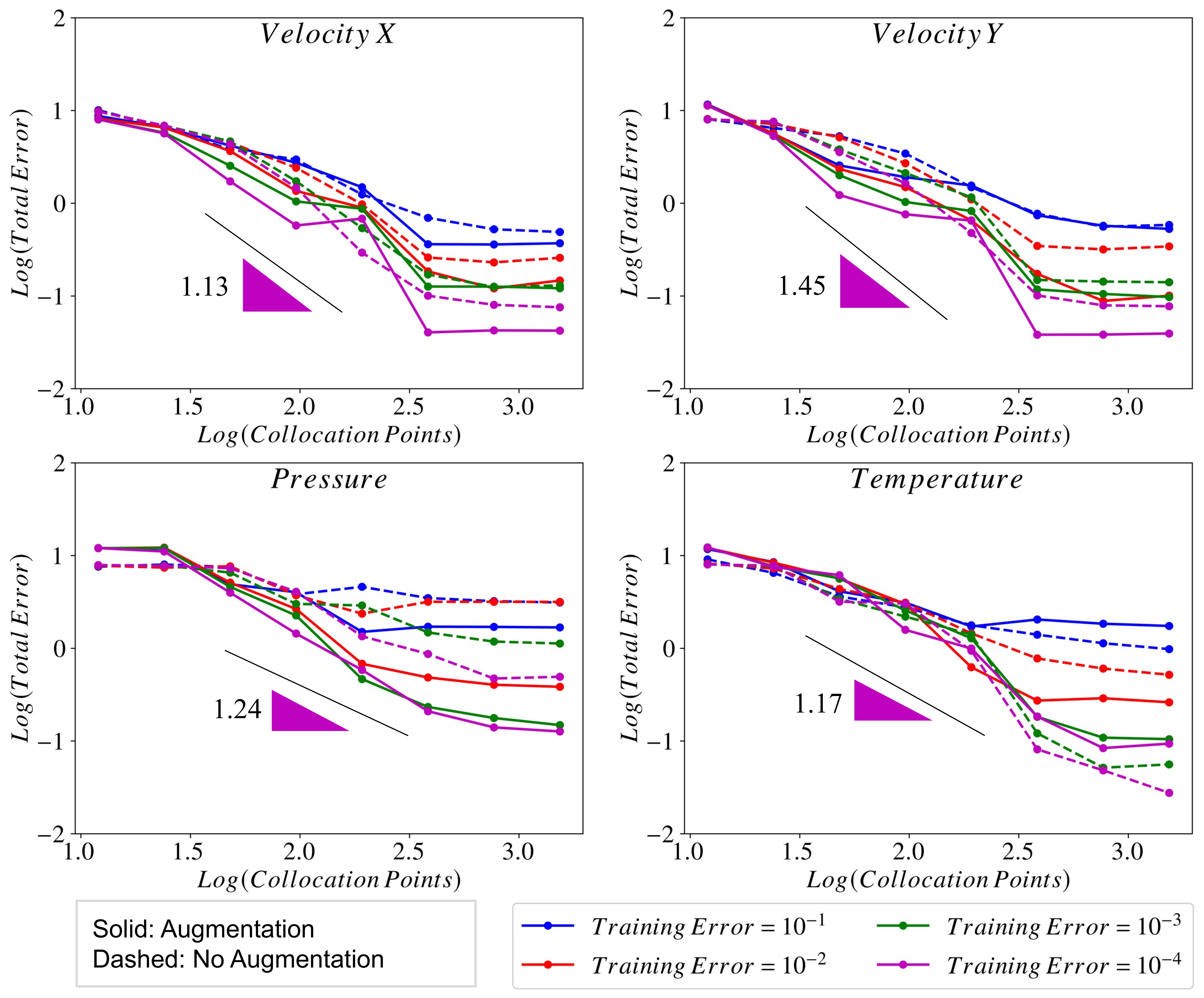}
    \caption{Convergence of \(W^{1,\infty}\) norm of error w.r.t training collocation points.}
    \label{fig: w1_conv}
    \end{figure}

\begin{figure}[H]
    \centering 
    \includegraphics[width= .7\linewidth]{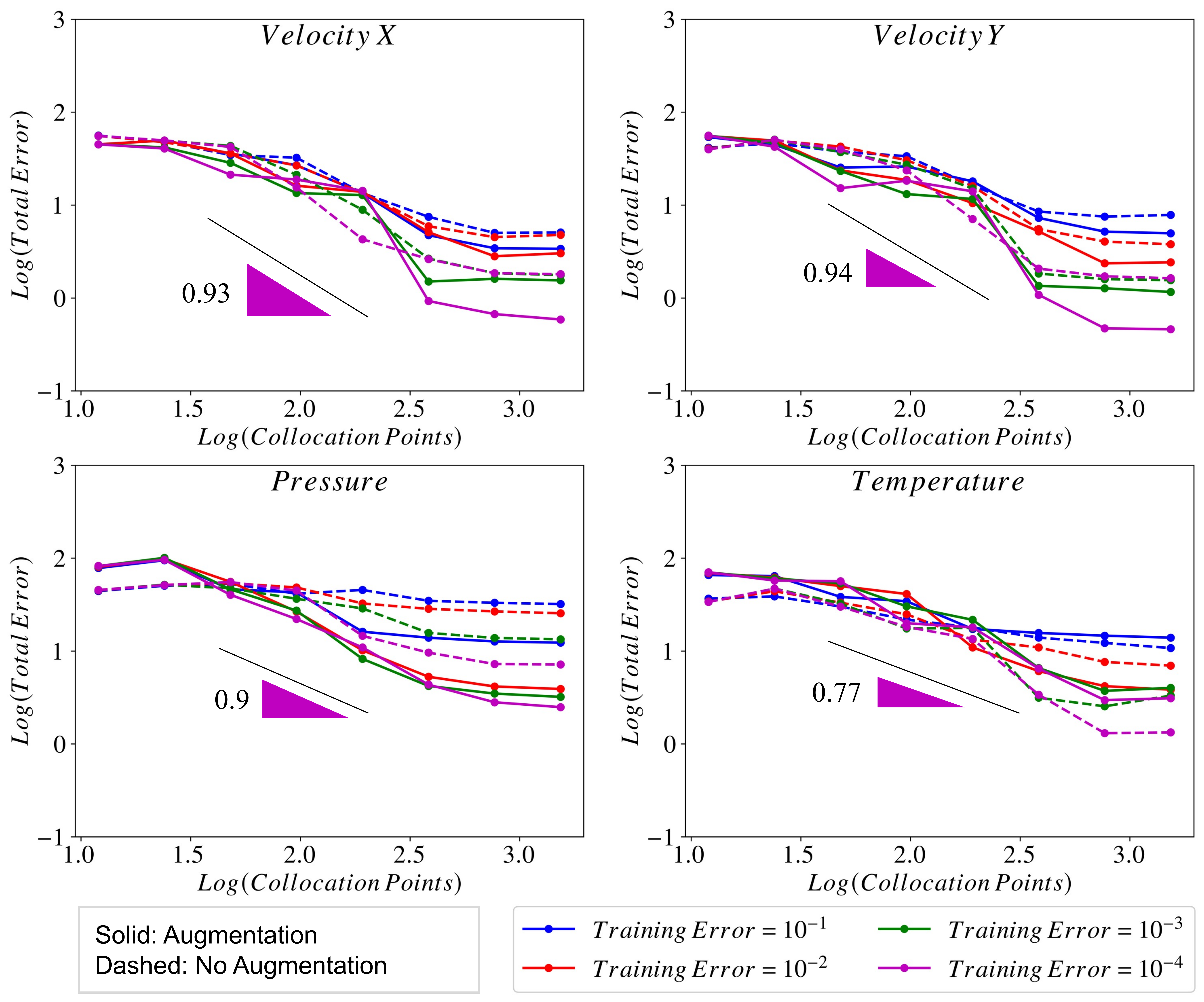}
    \caption{Convergence of \(W^{2,\infty}\) norm of error w.r.t training collocation points.}
    \label{fig: w2_conv}
    \end{figure}

    A consistent neural network converges to the exact solution if the sample size grows to infinity \cite{strang2019linear,kutz2013data}. This section explores the consistency of PINNs within a sample limit, which is equivalent to the convergence of the generalization error. As seen in the previous sections, number of collocation points affects the accuracy of PINNs predicted solution. Typically the number is selected using a trial-and-error procedure, however, a systematic approach is needed to sample enough number of points that can predict correct physics of the problem while minimizing the total error. For this purpose, an error convergence study w.r.t the training collocation points is carried out. The points are hierarchically increased such that the current set of points are fully represented in the subsequent larger training dataset, thereby leading to a consistent convergence analysis. Figures \ref{fig: w0_conv}-\ref{fig: w2_conv} show convergence plots for each of the velocity, pressure, and temperature fields and their gradients at various training residual thresholds. A good estimate for the number of collocation points is reached when the error convergence rate drops and there is not much reduction in error with the subsequent increase in collocation points.
    
    The convergence analysis shows that the rate of convergence increases for lower training error thresholds. The drop in the convergence rates indicate that the network has reached its approximation capacity with the given number of collocation points and the training error threshold. For the training error of \({10^{-4}}\), the convergence rates in the \(W^{0,\infty}\) norm of error for all the solution fields are super-linear i.e., between \(\mathcal{O}\left( {{N^{3/2}}}  \right)\) and \(\mathcal{O}\left( {{N^{2}}} \right)\). The convergence rates drop to \(\mathcal{O}\left( {{N}}  \right)\) for \(W^{1,\infty}\) and \(W^{2,\infty}\) norms for error that also includes first and second order derivatives of the solutions fields, respectively.
    
            
\begin{figure}[H]
    \centering 
    \includegraphics[width= .75\linewidth]{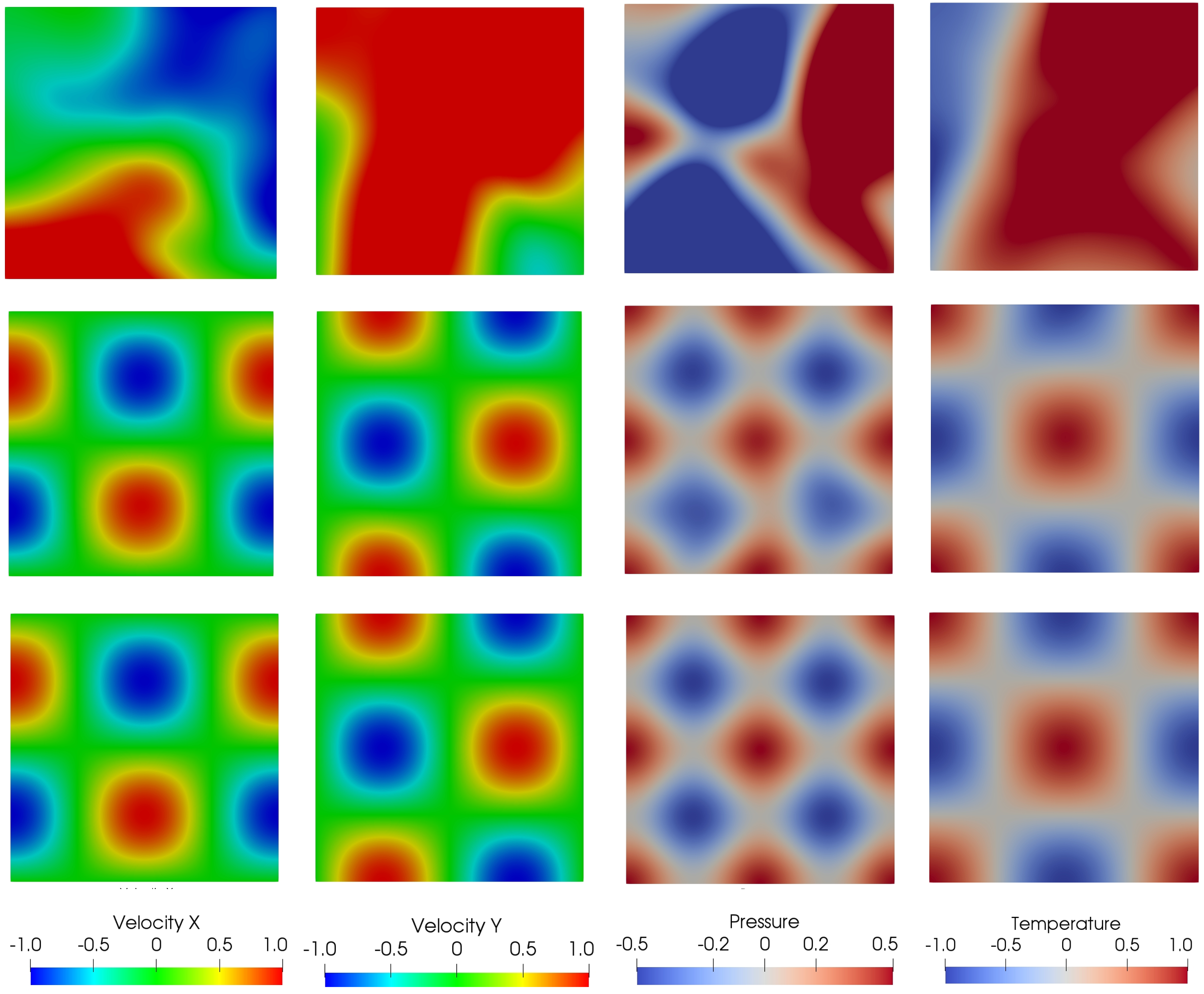}
    \caption{Qualitative plots of the solution fields for PINNs trained at different number of collocation points (first row = 12 points, second row = 96 points, third row = 1536 points).}
    \label{fig: conv_soln_qual}
    \end{figure}

\begin{figure}[H]
    \centering 
    \includegraphics[width= .75\linewidth]{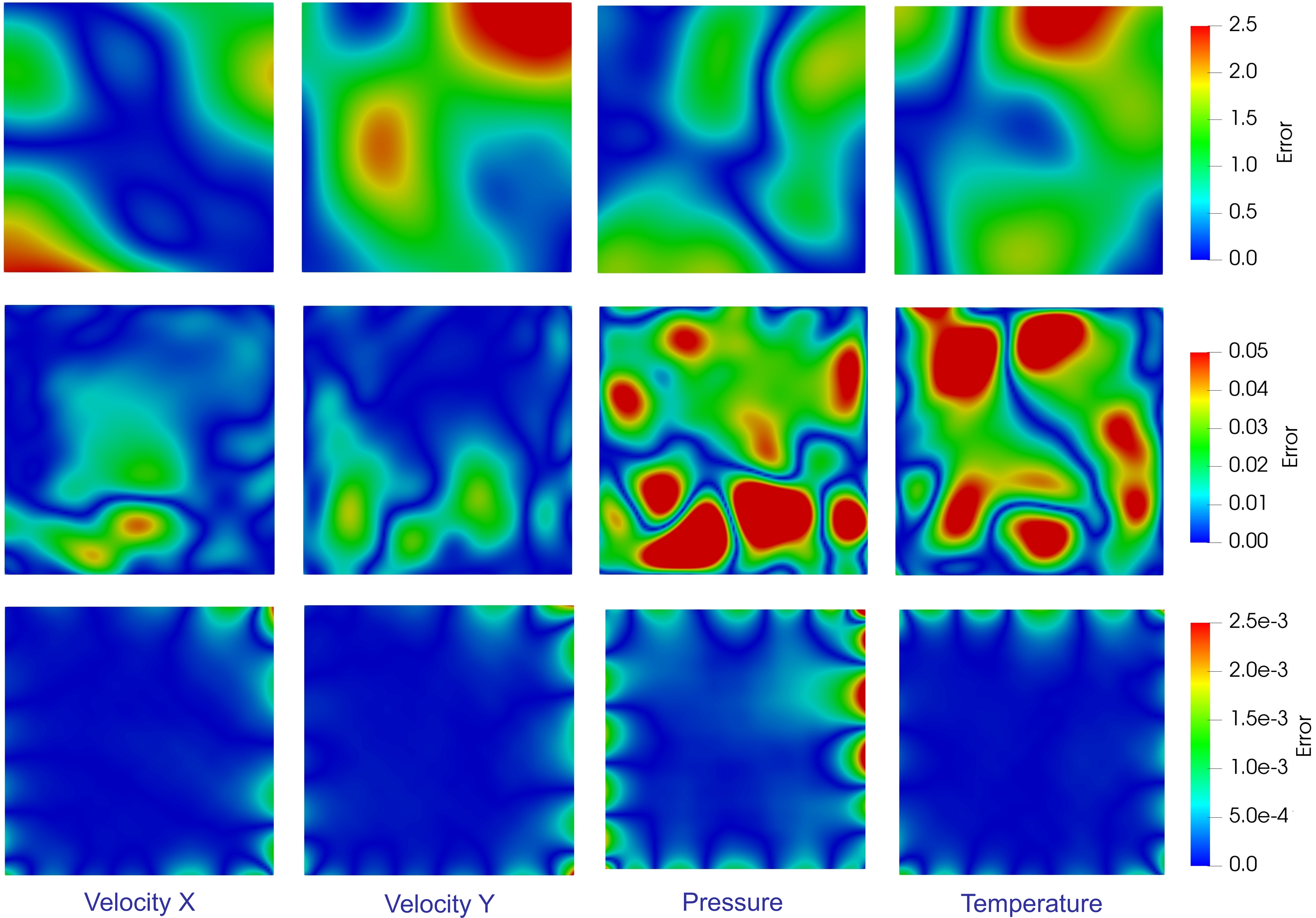}
    \caption{Qualitative plots of the pointwise error w.r.t the exact solution for PINNs trained at different number of collocation points (first row = 12 points, second row = 96 points, third row = 1536 points).}
    \label{fig: conv_err_qual}
    \end{figure}

    The PINNs predicted solution fields for various number of collocation points are shown in Figure \ref{fig: conv_soln_qual}. Although the network converges to the given training error at lower number of collocation points (first row in Figure \ref{fig: conv_soln_qual}), it fails to produce an accurate solution because it converges to wrong minima. Moreover, even though the qualitative solution fields in the second and third row look quite similar, the error plots in Figure \ref{fig: conv_err_qual} indicate that there is an order of magnitude difference in the error between the two solutions.  This highlights the importance of having a network that is consistent, and an error convergence study for selecting the appropriate number of collocation points for a particular problem. In Figure \ref{fig: conv_err_qual}, the error bar on the right indicates the magnitude of pointwise error in the predicted velocity, pressure, and temperature fields corresponding to that row.
    
    \begin{remark} \label{rem: 2}
    {Number of collocation points affects the accuracy of PINNs predicted solution and these points should be chosen depending on the length scale of the physics of the problem.}
    \end{remark}
    
\section{Generalization of Error Estimates for Different Problem Domain and Higher Reynolds Number Flow using Transfer Learning} \label{Sec: Generalization}

    In general, PINNs need to be trained again if flow parameters or problem domain is changed. Since the computational cost of training PINNs is very high, re-training from scratch is not cost effective for its generalization. However, the training of PINNs can be accelerated by using transfer learning. In transfer learning, instead of training from scratch, the parameters of the already trained network (e.g., at a lower Reynolds number) are transferred to a network for different flow conditions (e.g., a higher Reynolds number flow). In this way, only a fraction of the original epochs are needed to obtain an accurate solution. We have used \textit{BFGS} optimizer in transfer learning which is faster when the network parameters are already trained and they only need some fine-tuning.
    
\begin{figure}[H]
        \centering
        \begin{subfigure}[t]{.45\linewidth}
            \centering
            \includegraphics[width=\linewidth]{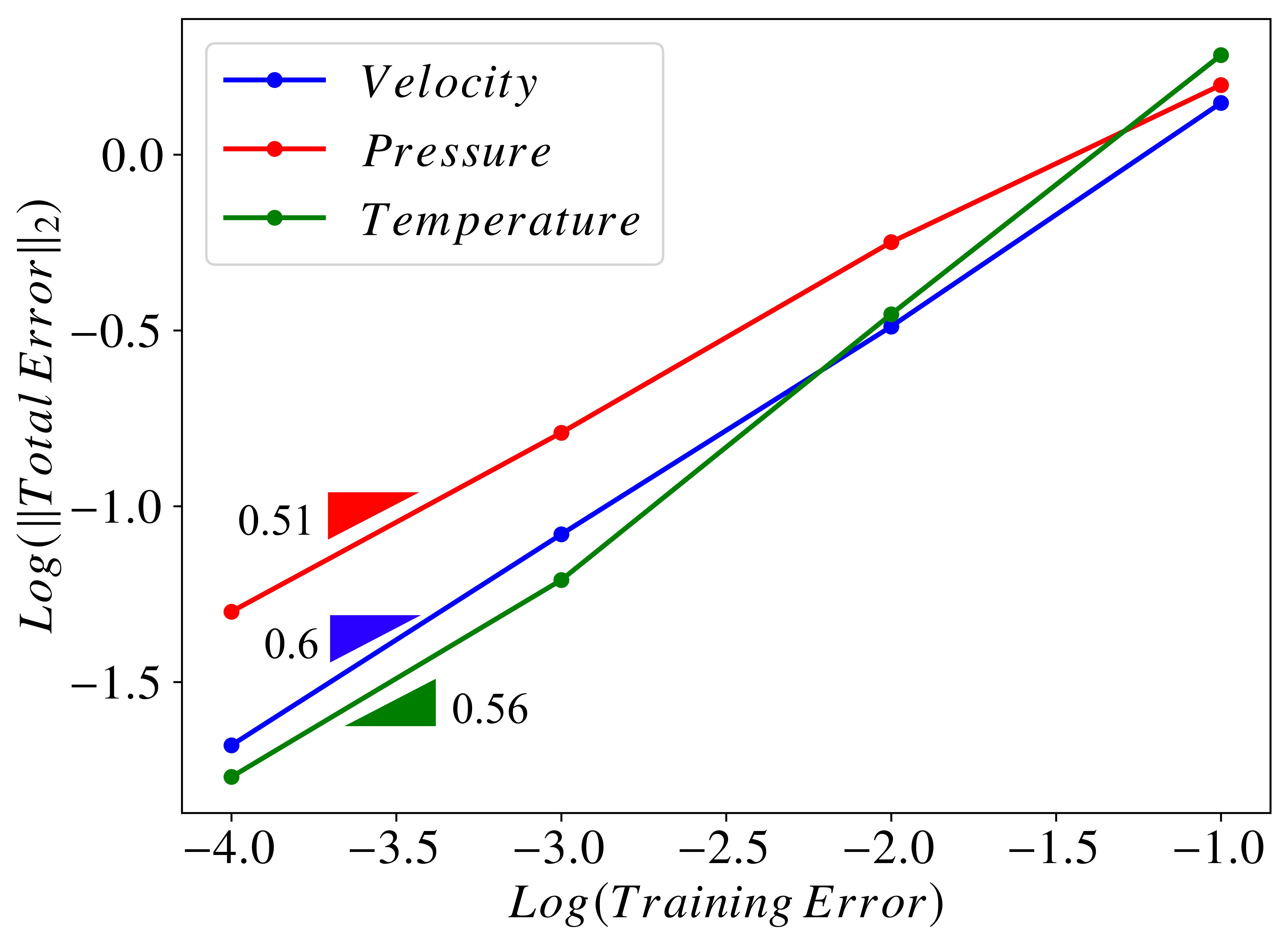}
            \caption{}
            \label{Fig10a}
        \end{subfigure}
        \hspace{1.5em}%
        \begin{subfigure}[t]{.45\linewidth}
            \centering
            \includegraphics[width=\linewidth]{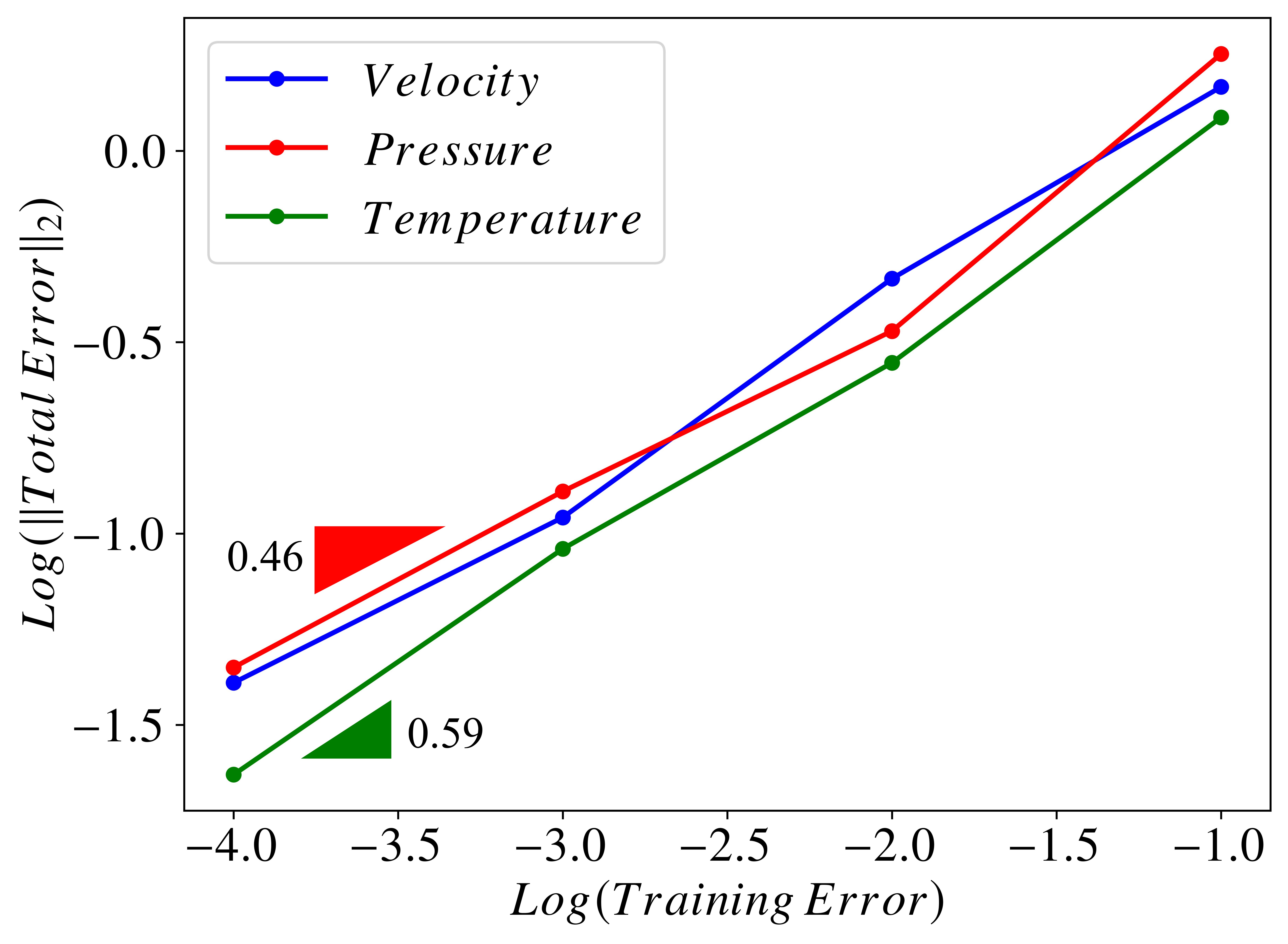}
            \caption{}
            \label{Fig10b}
        \end{subfigure}
        \caption{Scaling of \(L_2\) norm of total error with the training error for (a) different spatial domain of the problem (b) higher Reynolds number with a larger body force.}
        \label{fig: l2_geom_mat}
    \end{figure}

    To establish the generalization of the error estimates presented in Section \ref{Sec: NumericExp}, we employ transfer learning to train networks for two test cases. The first test case involves a different geometric description of the problem domain which is half of the original domain in the x-direction i.e., \(\Omega  = \left[ {0,1} \right] \times \left[ {-1,1} \right]\). The sinusoidal solution is no longer axisymmetric for the new geometric layout. In the second test case, we retain the spatial domain as that of the original problem, while we increase the Reynolds number from 1 to 10 and modify the value of gravitational acceleration from -1 to -9.8. In both the cases, the \(L_2\)-norm of the total error scales as the square root of the training error in accordance  with Eq. \ref{eq:29} as shown in Figure \ref{fig: l2_geom_mat}. Moreover, the convergence rates of total error w.r.t the training error and the number of collocation points are in the same ballpark as that for the original problem. These convergence rates are presented in \ref{App: A} and \ref{App: B}. These results show the practical significance of the numerical error estimates obtained in this work as they are consistent and can be generalized to different spatial domains or flow parameters in the coupled thermal flow problems. 

\section{Comparison between FEM and PINNs} \label{SubSec: CompFEMandPINN}
    In this section, we compare PINNs and a stabilized finite element method (FEM) \cite{zhu2022nano} for the solution of Beltrami flow. The FEM solution is computed on a \(100 \times 100\) structured mesh of linear quadrilateral and triangular elements with \(39,204\) degrees of freedom. For PINNs, 2-128-128-4 network is employed with 1536 total collocation points and a training error of \({10^{-4}}\). The trained network is then used to predict solution over a \(100 \times 100\) grid of testing points as shown in Figure \ref{fig: train_test_pts} and point-wise error w.r.t analytical solution is computed.
    
    \begin{figure}[H]
    \centering 
    \includegraphics[width= .8\linewidth]{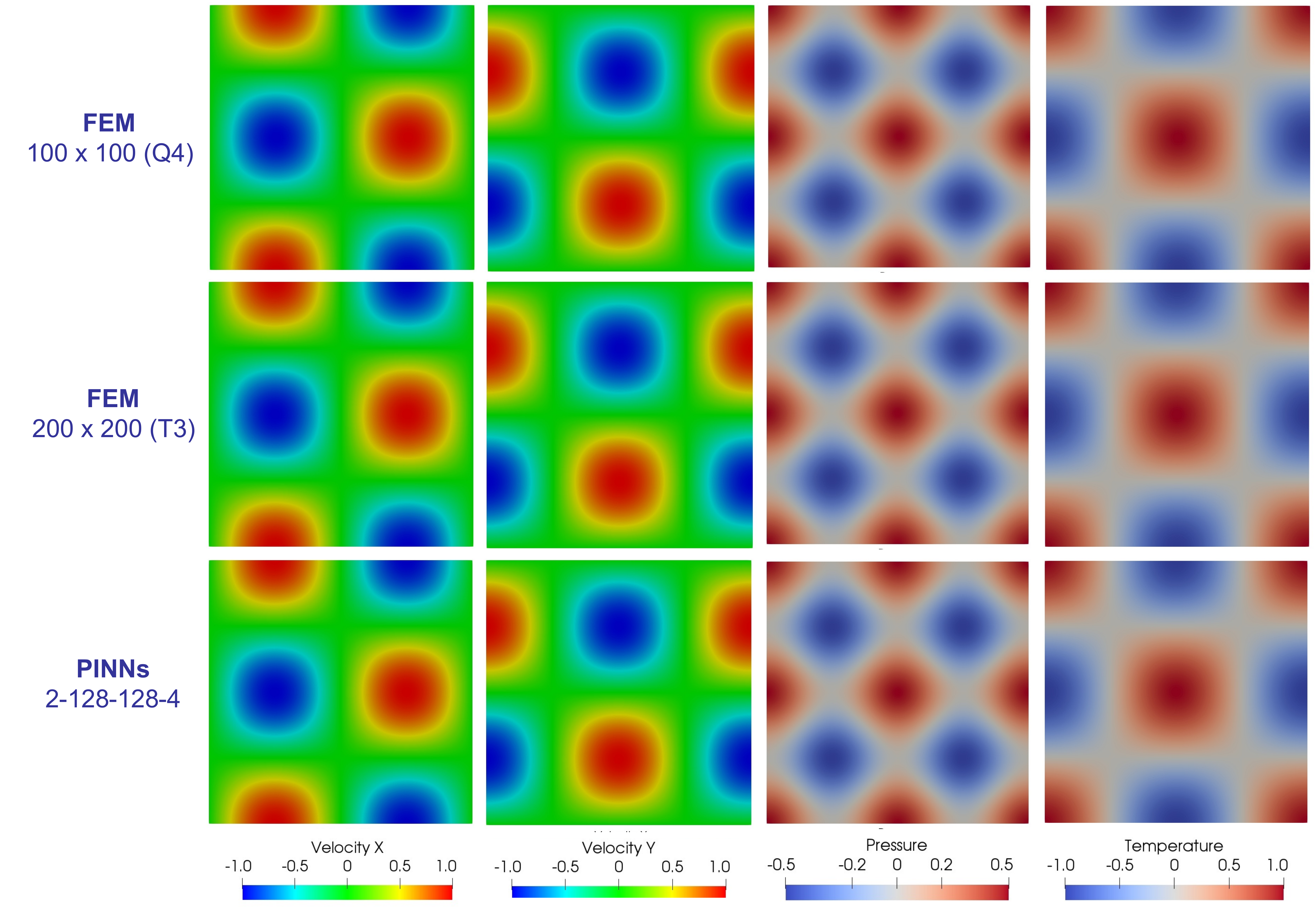}
    \caption{Comparison of solution fields for FEM and PINNs  (First row: FEM with \(100 \times 100\) linear quadrilateral elements; Second row: FEM with \(200 \times 200\) linear triangular elements; Third row: PINNs with 2-128-128-4 network and training error = \(10^{-4}\)).}
    \label{fig: comp_soln_FEM_PINN}
    \end{figure}
    
\begin{figure}[H]
    \centering 
    \includegraphics[width= .8\linewidth]{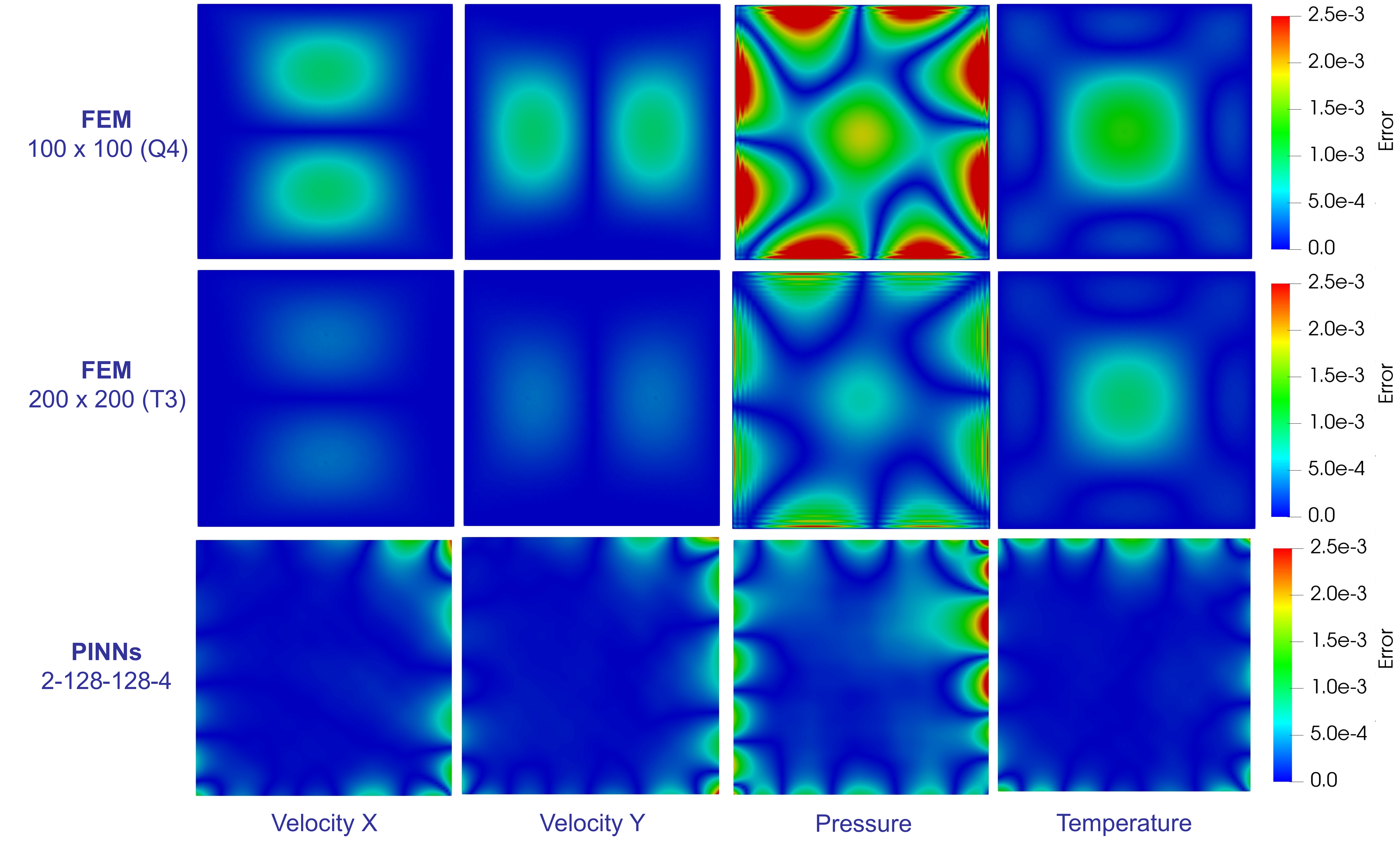}
    \caption{Comparison of error fields for FEM and PINNs (First row: FEM with \(100 \times 100\) linear quadrilateral elements; Second row: FEM with \(200 \times 200\) linear triangular elements; Third row: PINNs with 2-128-128-4 network and training error = \(10^{-4}\)) and collocation points = \(1536\).}
    \label{fig: comp_error_FEM_PINN}
    \end{figure}

    The solution fields as well as error w.r.t analytical solution for both the methods are presented in Figures \ref{fig: comp_soln_FEM_PINN} and \ref{fig: comp_error_FEM_PINN}, respectively. Qualitatively there is no significant difference between the solution fields from both the methods. However, PINNs show better error performance in the domain interior as the FEM solution has a larger magnitude of error in there as compared to PINNs. Particularly for the pressure field, PINNs predicted solution looks much better in the domain interior as compared to FEM which is due to the pressure augmentation term we have introduced in our model. On the boundaries, error from the FEM solution is lesser as compared to the error in PINNs. This is because the boundary conditions are strongly enforced in FEM while they are weakly enforced in PINNs. We also notice the non-symmetry of the error fields for the PINNs predicted solution which is due to the randomly chosen collocation points that are not symmetric. Selecting a structured grid of collocation points for training would result in symmetric error plots.
    
    Computationally PINNs take much more computational time and resources to train and predict the solution as compared to FEM. FEM takes about 25 seconds to generate the solution on a single CPU (excluding the time taken for meshing) while PINNs take almost 20 hours to train on a single CPU for the largest training dataset (1536 points) and a training residual of \(10^-4\). Once the network is trained, it predicts the solution instantaneously on any given set of points. Moreover, for a lower training error threshold and lesser number number of collocation points (e.g., \(10^{-3}\) and \(768\)), PINNs take only about an hour to train and the solution is good enough for all practical design purposes. The training time for PINNs can be further reduced by using GPUs, faster optimization algorithm such as \textit{BFGS}, and transfer learning. Although PINNs are computationally less efficient for forward problems with structured meshes, they can be efficient for problems where meshing takes considerable time.
        
\section{Conclusions}\label{Sec: Conclusion}
    
    We have presented convergence analysis for PINNs, applied to the solution of 2D steady incompressible Navier-Stokes PDEs coupled to the scalar energy PDE and subjected to Dirichlet boundary conditions. A pressure stabilization term in the form of pressure Poisson equation is added to the PDE residuals for PINNs. This physics informed augmentation of thermally coupled Navier-Stokes equations improves accuracy of the pressure field by an order of magnitude as compared to the case without augmentation. Via a model problem of Beltrami flows for which an analytical solution exists, we investigated the effects of (i) network size, (ii) collocation points, (iii) training error, and (iv) physics-informed augmentation of the model, on the reduction in the total error in the predicted solution. Through a carefully designed set of numerical test cases, \textit{posteriori} error estimates were obtained. The network convergence study showed better approximation capability and faster convergence of deeper networks that have larger trainable parameters, an outcome which is in accordance with the universal approximation theorem. The convergence of PINNs was established by showing that with enough training parameters and collocation points, a small training error leads to a small total error. The convergence rate of total error w.r.t the training error was observed to be sub-linear for different Sobolev norms, while the total error was observed to scale as the square root of training error in the \(L_2\) norm. The consistency of PINNs was shown through convergence study of total error w.r.t the training dataset. The rates of convergence of total error w.r.t the number of collocation points were super-linear for all the solution fields (\(W^{0,\infty}\) norm) and were linear for their first and second derivatives (\(W^{1,\infty}\) and \(W^{2,\infty}\) norms). Faster convergence rates were observed as the training error threshold was reduced. Moreover, we showed that all these convergence rates can be generalized for problems with different spatial domain, and/or different flow parameters such as a larger body force and a higher Reynolds number within the laminar range. A comparison of the PINNs predicted solution with the finite element based solution also highlighted the various mathematical attributes of PINNs for application to coupled multiphysics problems. 
\section*{Acknowledgments}
    Computing resources for the numerical test cases were provided by the Teragrid/XSEDE Program under NSF Grant TG-DMS100004. 

\bibliographystyle{elsarticle-num}
\bibliography{references}

\appendix
\section{Error Estimates for Different Problem Domain} \label{App: A}

    \setcounter{figure}{0}

    \begin{figure}[H]
    \centering \includegraphics[width=.7\linewidth]{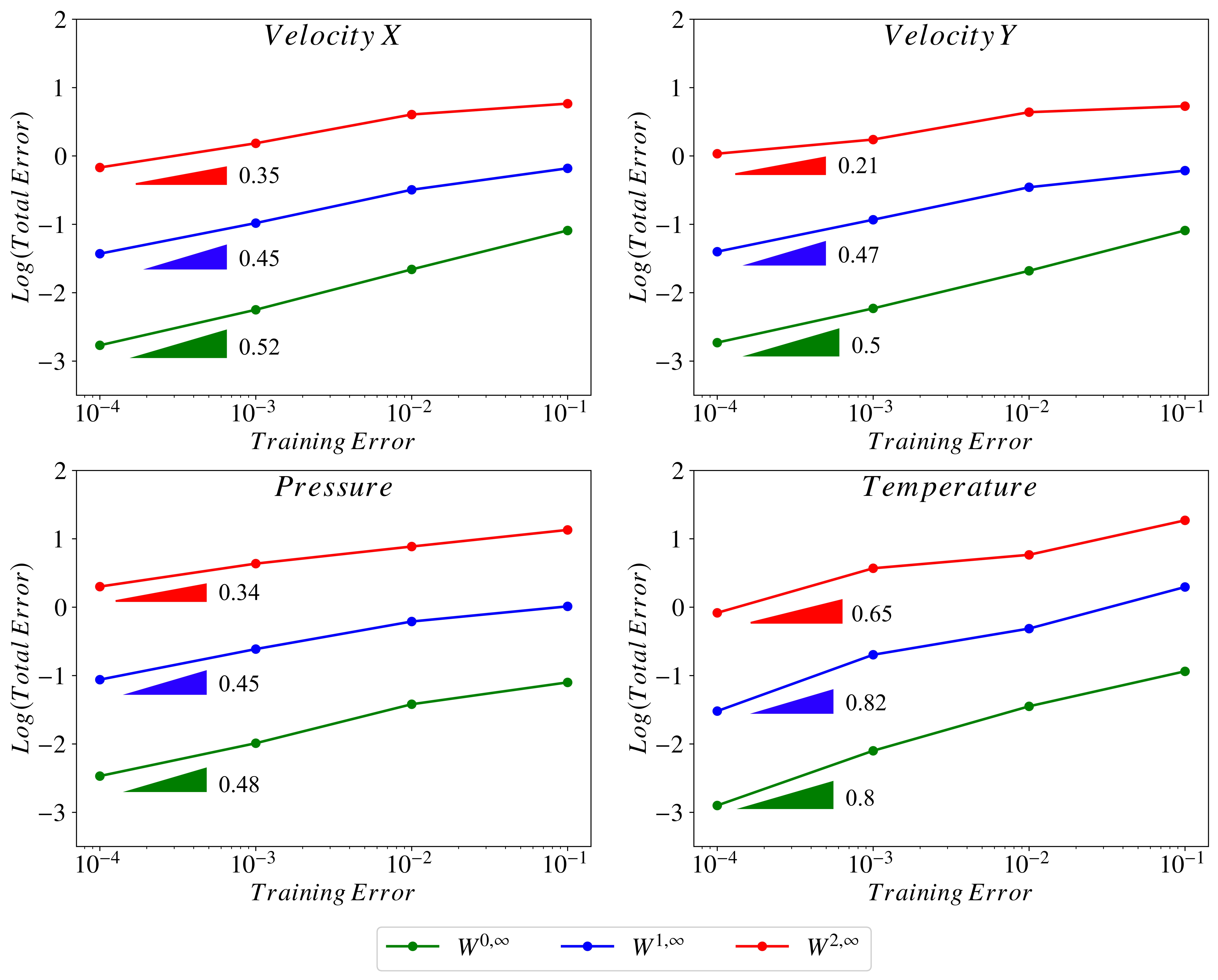}
    \caption{\small Convergence of total error w.r.t training error under different Sobolev norms.}
    \label{fig: geom_train}
    \end{figure}
    
    \begin{figure}[H]
    \centering \includegraphics[width=.7\linewidth]{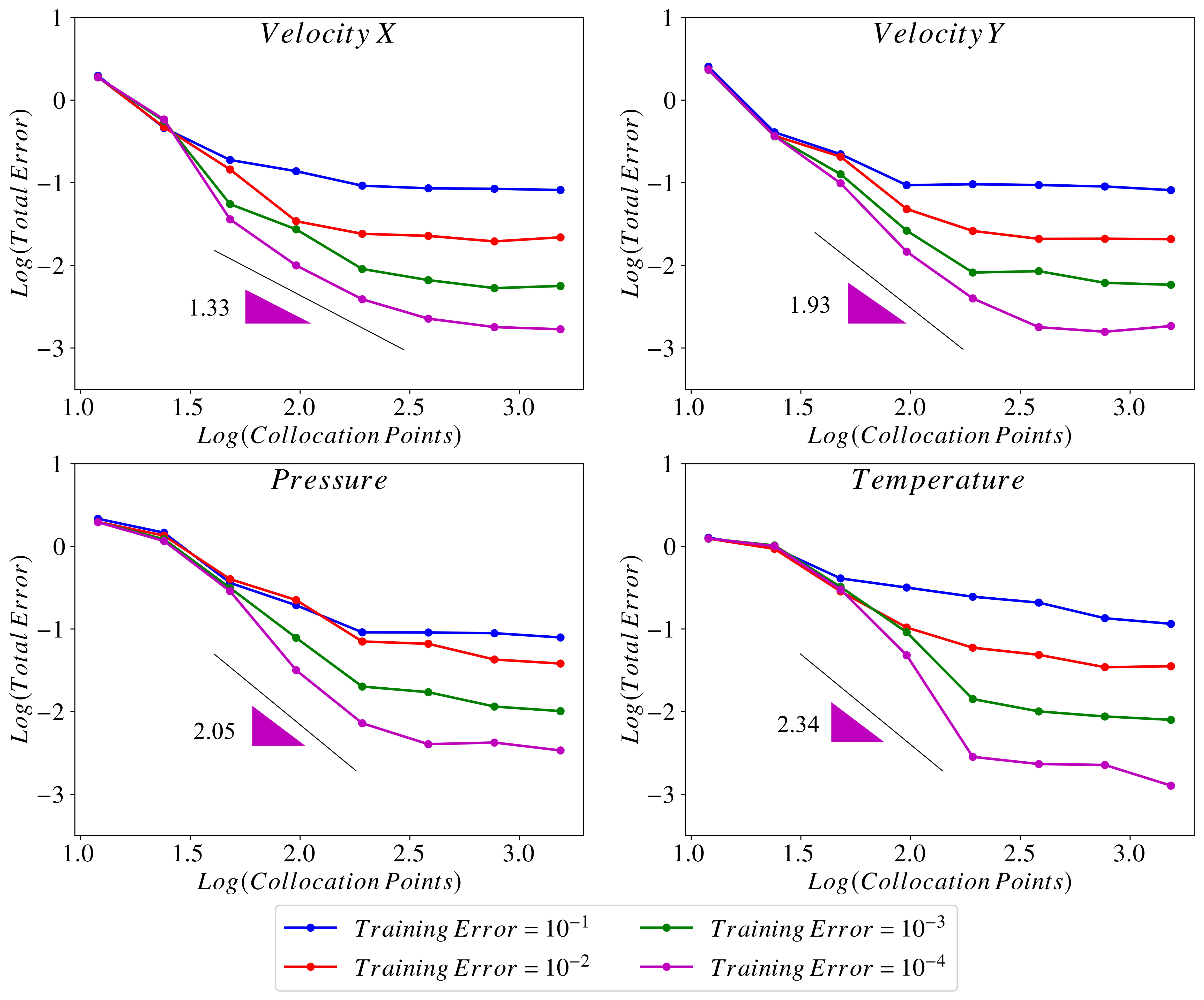}
    \caption{\small Convergence of \(W^{0,\infty}\) norm of error w.r.t the number of training collocation points.}
    \label{fig: geom_w0}
    \end{figure}
    
    \begin{figure}[H]
    \centering \includegraphics[width=.7\linewidth]{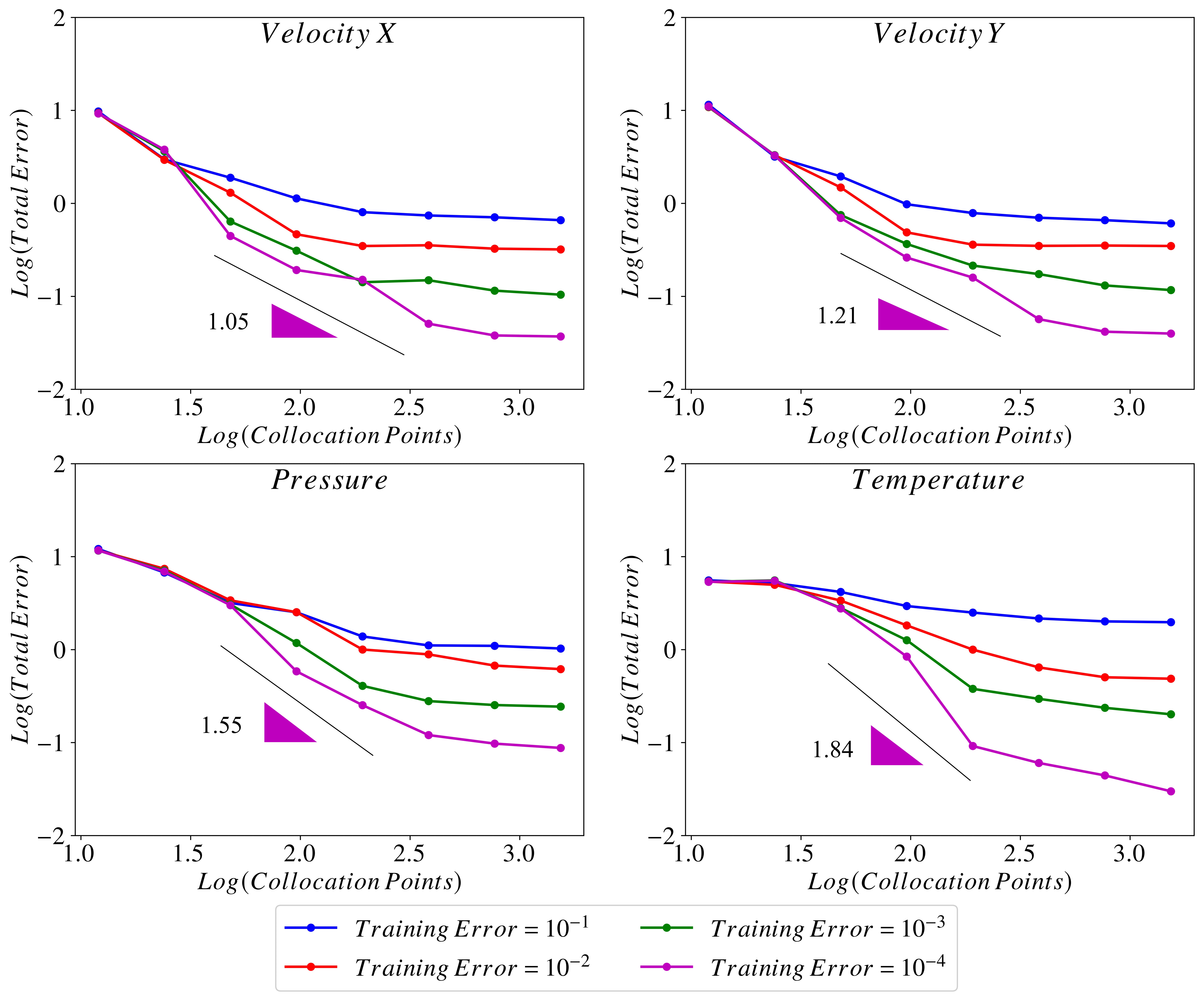}
    \caption{\small Convergence of \(W^{1,\infty}\) norm of error w.r.t the number of training collocation points.}
    \label{fig: geom_w1}
    \end{figure}
    
    \begin{figure}[H]
    \centering \includegraphics[width=.7\linewidth]{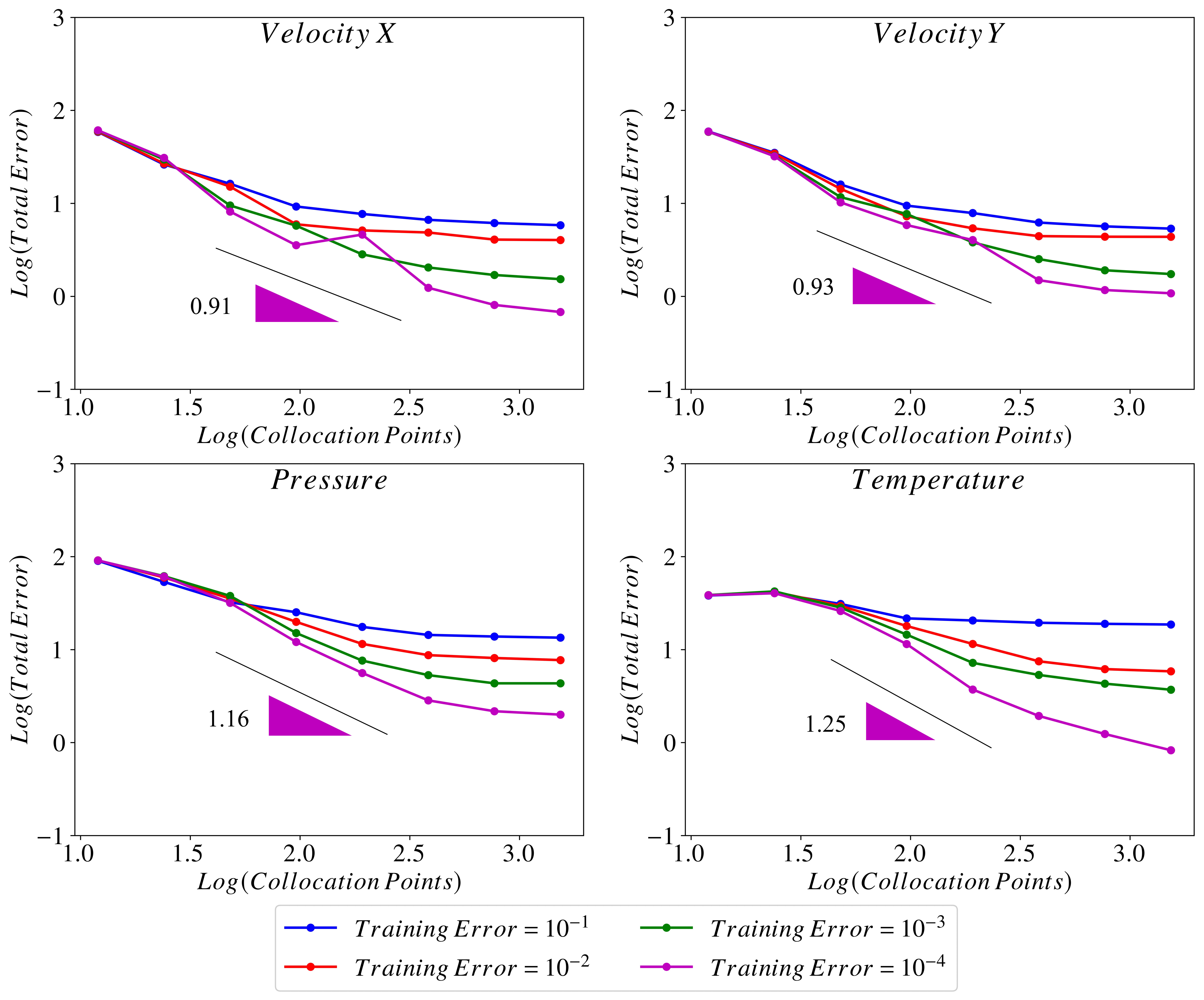}
    \caption{\small Convergence of \(W^{2,\infty}\) norm of error w.r.t the number of training collocation points.}
    \label{fig: geom_w2}
    \end{figure}

\section{Error Estimates for Higher Reynolds Number} \label{App: B}

    \setcounter{figure}{0}

    \begin{figure}[H]
    \centering \includegraphics[width=.7\linewidth]{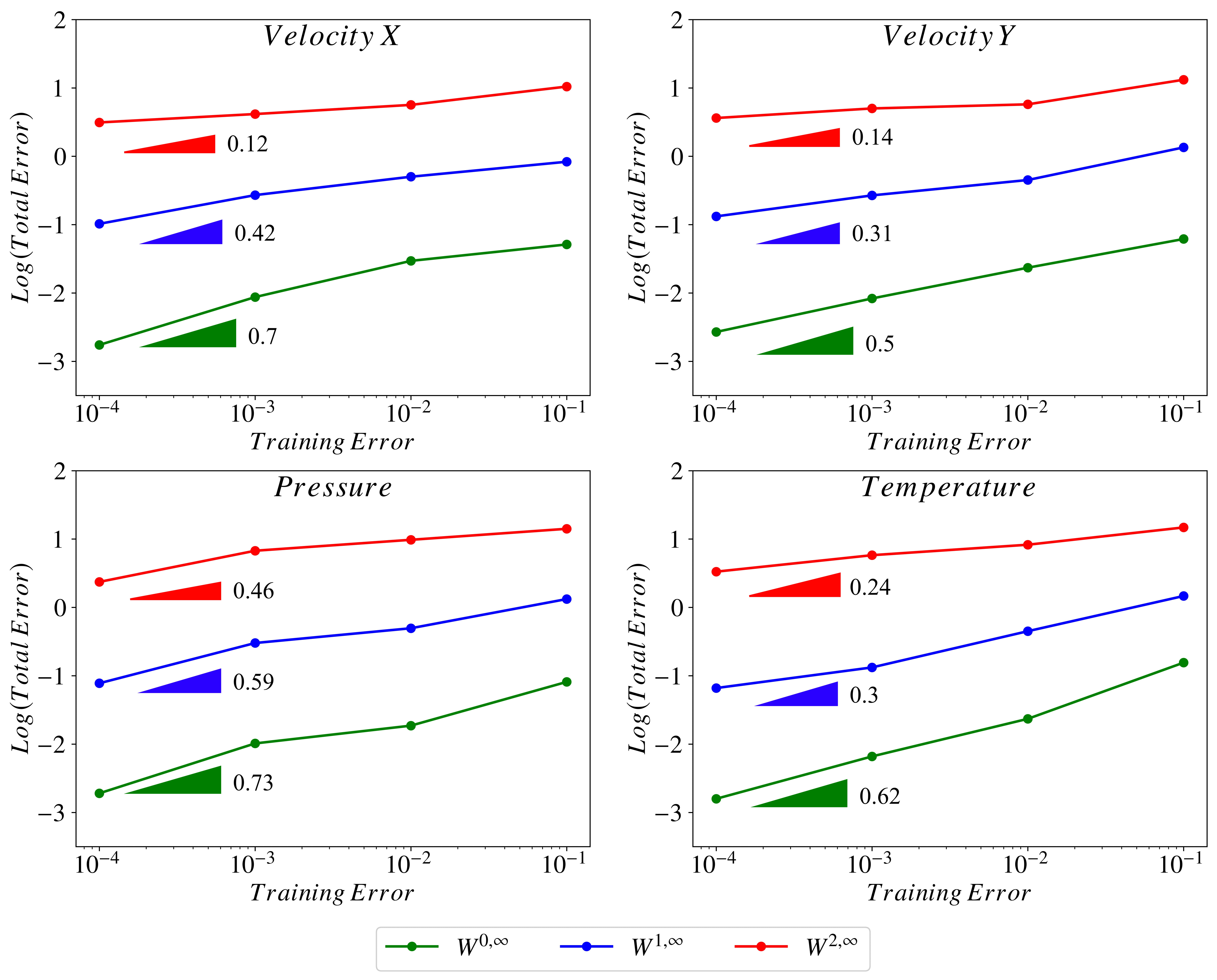}
    \caption{\small Convergence of total error w.r.t training error under different Sobolev norms.}
    \label{fig: mat_train}
    \end{figure}
    
    \begin{figure}[H]
    \centering \includegraphics[width=.7\linewidth]{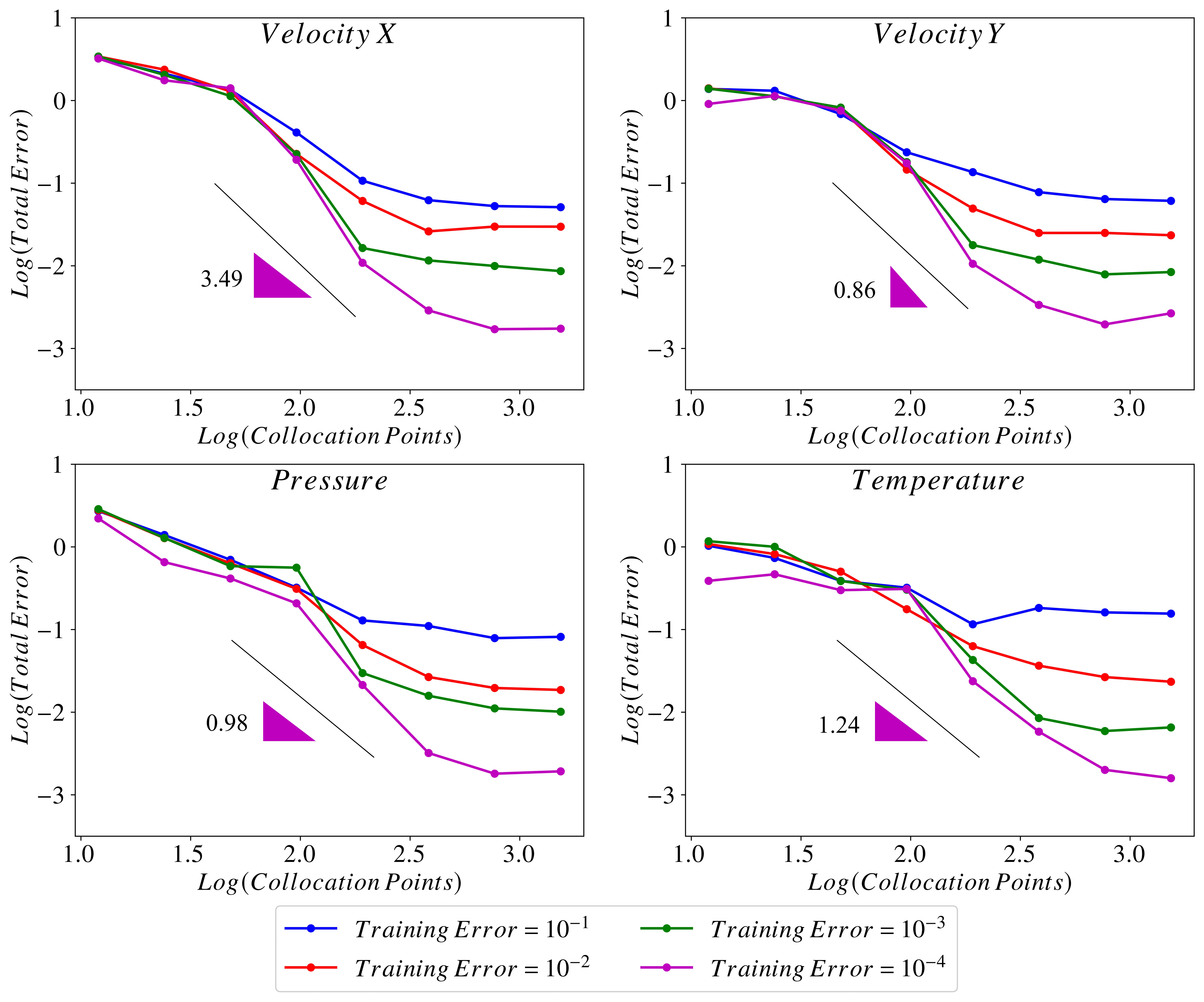}
    \caption{\small Convergence of \(W^{0,\infty}\) norm of error w.r.t the number of training collocation points.}
    \label{fig: mat_w0}
    \end{figure}
    
    \begin{figure}[H]
    \centering \includegraphics[width=.7\linewidth]{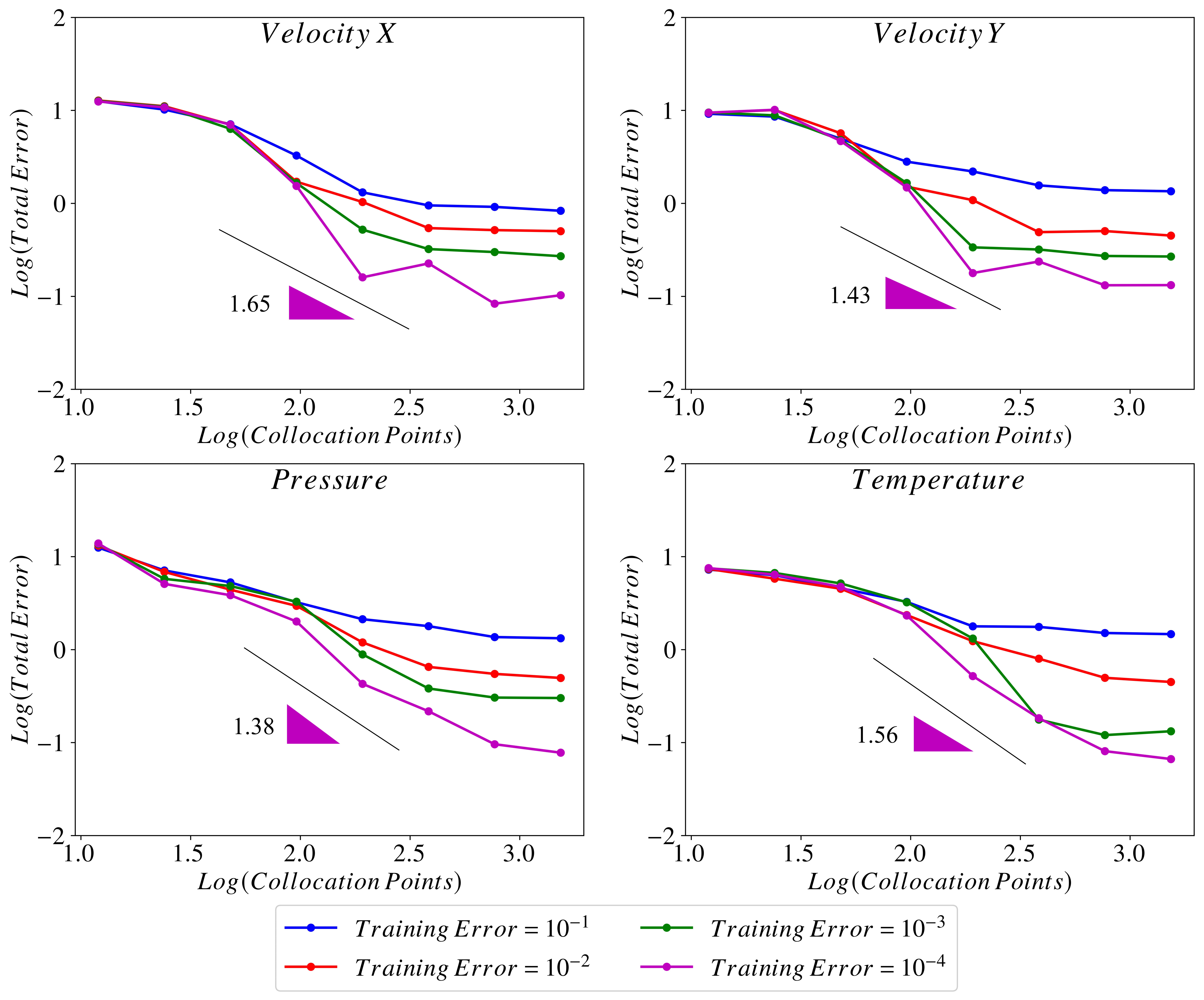}
    \caption{\small Convergence of \(W^{1,\infty}\) norm of error w.r.t the number of training collocation points.}
    \label{fig: mat_w1}
    \end{figure}
    
    \begin{figure}[H]
    \centering \includegraphics[width=.7\linewidth]{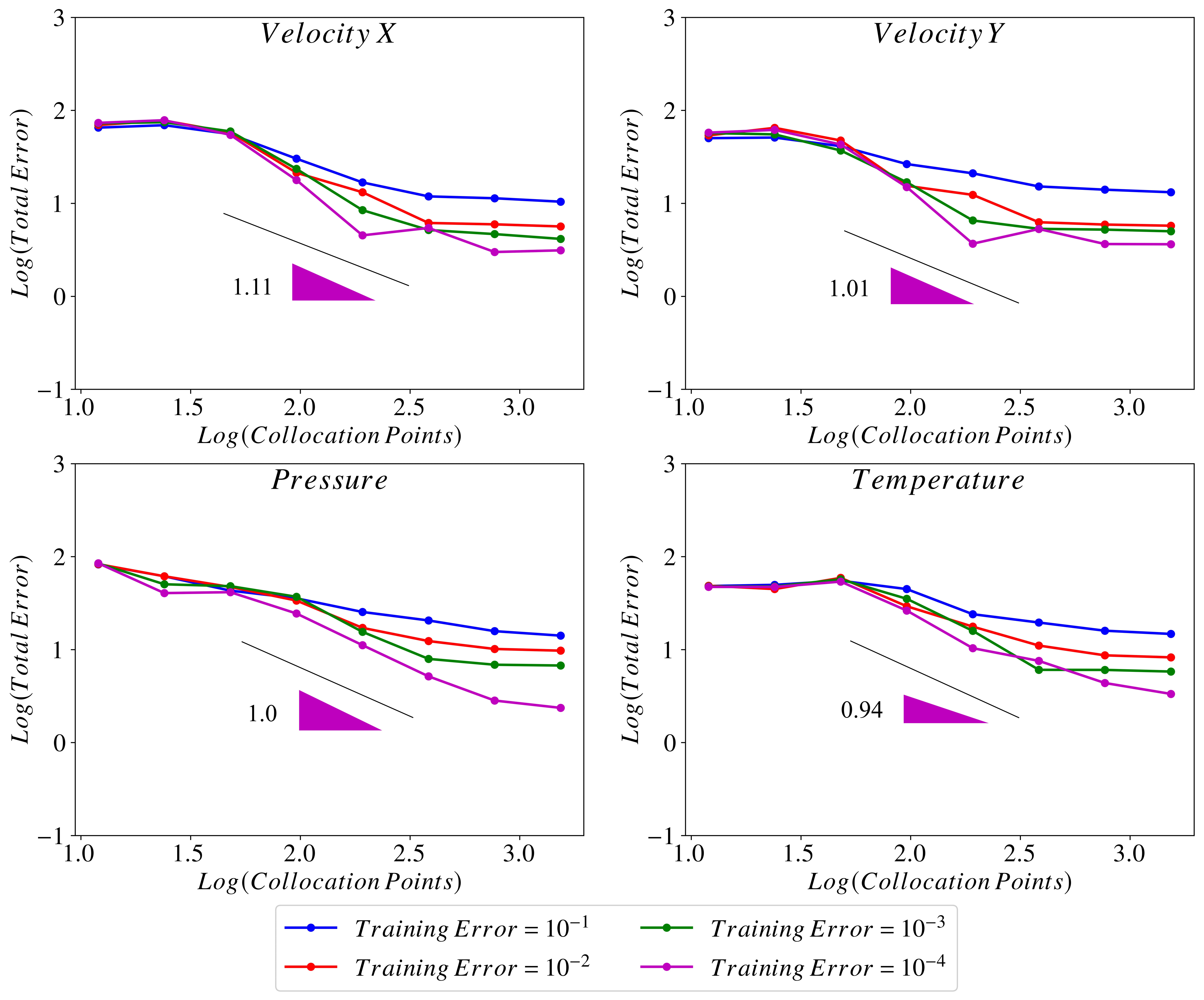}
    \caption{\small Convergence of \(W^{2,\infty}\) norm of error w.r.t the number of training collocation points.}
    \label{fig: mat_w2}
    \end{figure}

\end{document}